\documentclass[10pt, reqno, hidelinks]{amsart}
\usepackage{amsthm,amssymb,amscd,amsmath,amsfonts,amssymb,amscd,stmaryrd,hyperref,mathrsfs,yhmath,mathtools}                                                                                                                    
\usepackage[all]{xy}
\usepackage{hyperref}
\usepackage{todonotes}
\usepackage{enumerate}
\usepackage{physics}
\usepackage{ragged2e}
\usepackage{xfrac, faktor} 
\usepackage{lipsum}  
\usepackage{tikz-cd}
\usepackage{imakeidx}
\usepackage{subfiles}
\usepackage{bbm}
\usepackage{quiver}
\usepackage{mdframed}
\usepackage[margin = 1.3in]{geometry}
\usepackage{changepage}

\usepackage{amsmath,amsfonts,amscd,amssymb,verbatim,amsthm}

\newcommand{\OO}{\mathcal{O}}

\DeclareMathOperator{\Isoc}{Isoc}

\DeclareMathOperator{\ind}{ind}

\DeclareMathOperator{\opp}{op}

\DeclareMathOperator{\Gal}{Gal}

\DeclareMathOperator{\Br}{Br}

\DeclareMathOperator{\charfield}{char}

\DeclareMathOperator{\Vect}{\bold{Vect}}

\DeclareMathOperator{\Rep}{Rep}
\DeclareMathOperator{\End}{End}

\DeclareMathOperator{\inv}{inv}
\DeclareMathOperator{\Frob}{Frob}

\DeclareMathOperator{\ord}{ord}

\DeclareMathOperator{\Ind}{Ind}
\DeclareMathOperator{\Irr}{Irr}
\DeclareMathOperator{\GL}{GL}
\DeclareMathOperator{\ab}{ab}

\DeclareMathOperator{\Aut}{Aut}

\DeclareMathOperator{\Spec}{Spec}

\newcommand{\twopartdef}[4]
{
	\left\{
		\begin{array}{ll}
			#1 & \mbox{if } #2 \\
			#3 & \mbox{if } #4
		\end{array}
	\right.
}

\makeatletter
\newcommand{\colim@}[2]{%
  \vtop{\m@th\ialign{##\cr
    \hfil$#1\operator@font colim$\hfil\cr
    \noalign{\nointerlineskip\kern1.5\ex@}#2\cr
    \noalign{\nointerlineskip\kern-\ex@}\cr}}%
}
\newcommand{\colim}{%
  \mathop{\mathpalette\colim@{\rightarrowfill@\textstyle}}\nmlimits@
}
\makeatother

\newcommand{\bC}{{\mathbb C}}

\newcommand{\bF}{{\mathbb F}}

\newcommand{\bH}{{\mathbb H}}

\newcommand{\bN}{{\mathbb N}}

\newcommand{\bQ}{{\mathbb Q}}
\newcommand{\bR}{{\mathbb R}}
\newcommand{\bZ}{{\mathbb Z}}

\newcommand{\sB}{{\mathcal B}}
\newcommand{\sC}{{\mathcal C}}
\newcommand{\sD}{{\mathcal D}}

\newcommand{\sN}{{\mathcal N}}

\newcommand{\sP}{{\mathcal P}}

\newcommand{\fC}{{\mathfrak C}}

\newcommand{\fm}{{\mathfrak m}}
\newcommand{\fp}{{\mathfrak p}}

\newtheorem{thm}{Theorem}
\numberwithin{thm}{subsection}
\newtheorem{lemma}[thm]{Lemma} 
\newtheorem{prop}[thm]{Proposition} 
\newtheorem{cor}[thm]{Corollary}

\newtheorem*{theoremA}{Theorem A}

\theoremstyle{definition}
\newtheorem{defn}[thm]{Definition}
\newtheorem*{defn*}{Definition}
\newtheorem{eg}[thm]{Example}
\newtheorem*{eg*}{Example}
\newtheorem{remark}[thm]{Remark}

\title[A Normal Form for Invertible Semilinear Maps]{A Normal Form for Invertible Semilinear Maps and the Category of Isocrystals Over a Finite Field}
\date{\today}
\author{Yuewen Qi}
\address{Centre for Mathematical Sciences, University of Cambridge, Cambridge, UK}
\email{yq284@cam.ac.uk}
\author{James Taylor}
\address{School of Mathematics, University of Bristol, Bristol, UK}
\email{james.h.a.taylor@bristol.ac.uk}

\begin{document}
\maketitle

\begin{abstract}
For a field $L$, and a finite order automorphism $\sigma$, we classify invertible $\sigma$-semilinear endomorphisms on finite-dimensional $L$-vector spaces up to equivalence. This is deduced as a consequence of a description, in terms of certain norm equations in finite extensions of $L$, of the semilinear Schur indices for the group $\bZ$ acting on $L$ via $\sigma$. As an application, we give an explicit description of the category of isocrystals over a finite field.
\end{abstract}

\vspace{-0.3cm}

\setcounter{tocdepth}{1}
\tableofcontents

\vspace{-1.2cm}

\section{Introduction}

Suppose that $G$ is a group, $L$ is a field, and $G$ acts on $L$, encoded in a group homomorphism $\sigma \colon G \rightarrow \Aut(L)$ to the group of field automorphisms of $L$.

In this setting, one can consider the category $\Rep_L^{\sigma}(G)$ of semilinear representations of $G$ over $L$, whose objects are finite-dimensional $L$-vector spaces $W$, with an action $\rho$ of $G$, where now as opposed to linear representations, each $\rho(g)$ is required to be a $\sigma(g)$-linear endomorphism of $W$.

\subsection*{Motivation}

Semilinear representations and semilinear constructions arise naturally in many areas of mathematics, and such semilinear representation categories are worthy of systematic study.

For example, in $p$-adic Hodge theory and its precursors, the fundamental objects are semilinear: Dieudonn\'{e} theory classifies $p$-divisible groups over a perfect field $k$ of characteristic $p$ by finite free $W(k)$-modules equipped with a Frobenius-semilinear operator, and Fontaine's theory of $(\varphi,\Gamma)$-modules describes the category of all $p$-adic Galois representations of $G_K$, for $K/\mathbb{Q}_p$ finite, in terms of semilinear representations of $\varphi^{\bN} \times \Gamma$ over a certain ring $\mathbf{A}_K$ attached to $K$.

Semilinear representations also appear in algebraic geometry. For example, for a scheme $X$ with an action of a group $G$, it is a natural and fundamental problem to understand the category $\Vect^G(X)$ of $G$-equivariant vector bundles on $X$. In the most basic case that $X = \Spec(L)$ is a point, taking global sections canonically identifies $\Vect^G(X)$ with $\Rep_L^{\sigma}(G)$, and thus $\Rep_L^{\sigma}(G)$ provides a key example for how these categories behave in general.

\subsection*{Description in Terms of Linear Representations} Let $H \coloneqq \ker(\sigma \colon G \rightarrow \Aut(L)) \triangleleft G$ and $K \coloneqq L^G$.
In recent work of the first author \cite{TAYSchur}, it was shown how $\Rep_L^{\sigma}(G)$ can be described in terms of the linear representation category $\Rep_L(H)$, under the assumption that $L / K$ is finite. This assumption is automatic in many cases, such as if $G$ is finite, or $L$ is a number field.

We may summarise the main points as follows:

\begin{adjustwidth}{-0.0cm}{0cm}
\begin{itemize}
	\item Every object of $\Rep_L^{\sigma}(G)$ is uniquely the direct sum of indecomposable objects, so it suffices to describe $\Ind_L^{\sigma}(G)$, the set of isomorphism classes of indecomposable objects of $\Rep_L^{\sigma}(G)$.
	\item Towards describing this set, there is a natural action of $\Gal(L/K)$ on $\Ind_L(H)$, and for each $V \in \Ind_L(H)$ a number $m_L^{\sigma}(V) \in \bZ_{\geq 1}$, the semilinear Schur index of $V$.
	\item  These $m_L^{\sigma}(V)$ are constant on the $\Gal(L/K)$-orbits of $\Ind_L(H)$, and for any such orbit $\OO$ we can define $m_L^{\sigma}(\OO)$ to be $m_L^{\sigma}(V)$ for any representative $V \in \OO$.
	\item The main result of \cite{TAYSchur} is that $\Ind_L^{\sigma}(G)$ is indexed by the $\Gal(L/K)$-orbits of $\Ind_L(H)$.
	\item Explicitly, for each orbit $\OO$, $m_L^{\sigma}(\OO)$-copies of the orbit sum $\oplus_{V \in \OO} V$ extends uniquely (up to isomorphism) to a semilinear representation $W_{\OO}$ of $G$ over $L$, $W_{\OO}$ is indecomposable, and this $W_{\OO} \in \Ind_L^{\sigma}(G)$ is the corresponding indecomposable semilinear representation of $G$.
\end{itemize}
\end{adjustwidth}

Therefore, in practice, describing the category $\Rep_L^{\sigma}(G)$ can be broken down into two parts:
\begin{enumerate}
	\item Describe the $\Gal(L/K)$-set $\Ind_L(H)$,
	\item Determine the semilinear Schur indices $m_L^{\sigma}(V)$.
\end{enumerate}

The first is often readily computable, but the second is difficult in general. For example, if $G$ is finite and $L$ has characteristic zero, then the action of $\Gal(L/K)$ on $\Irr_L(H) = \Ind_L(H)$ can be computed from the character table, as some twisted variant of the standard Galois action, whilst the values of $m_L^{\sigma}(V)$ for each irreducible character $\chi_V$ are much more mysterious, and in general depend heavily on the arithmetic of $L$.

\subsection*{The Case When $G = \bZ$}
In this paper we consider the case when $G = \bZ$.

In this case, the data of the action map $\sigma \colon \bZ \rightarrow \Aut(L)$ corresponds to choosing the image of $1 \in \bZ$: from now on let $\sigma$ denote this automorphism of $L$, and write $\Rep_L^{\sigma}(\bZ)$ for the associated semilinear representation category. The assumption that $L /K$ is finite translates to the property that $\sigma \in \Aut(L)$ has finite-order $r \geq 1$, and $H$ is now simply the subgroup $r \bZ$.

From the results outlined above, to describe $\Rep_L^{\sigma}(\bZ)$, it suffices to
\begin{enumerate}
	\item Describe the $\Gal(L/K)$-set $\Ind_L(r\bZ)$,
	\item Determine the semilinear Schur indices $m_L^{\sigma}(V)$ for $V \in \Ind_L(r\bZ)$.
\end{enumerate}

For the first, there is a canonical identification
\[
	\Phi \colon \sP_L \times \bZ_{\geq 1} \xrightarrow{\sim} \Ind_L(r\bZ),
\]
where $\sP_L$ denotes the set of monic irreducible polynomials $p(x) \in L[x]$ with $p(x) \neq x$. Here the pair $(p(x),t)$ corresponds to the $\bZ$-representation
\[
V_{p(x), t, r} \coloneqq L[x] / p(x)^t,
\]
where $r \in \bZ$ acts by left multiplication by $x$. Furthermore, under the natural action of $\Gal(L/K)$ on $\sP_L$ through the coefficients and the trivial action on $\bZ_{\geq 1}$, $\Phi$ is $\Gal(L/K)$-equivariant.

We are therefore reduced to (2), determining $m_L^{\sigma}(V)$ for $V \in \Ind_L(r\bZ)$, or equivalently
\[
m_L^{\sigma}(p(x)^t) \coloneqq m_L^{\sigma}(V_{p(x), t, r}),
\]
for all $(p(x),t) \in \sP_L \times \bZ_{\geq 1}$ (after indexing $\Ind_L(r\bZ)$ by $\sP_L \times \bZ_{\geq 1}$).

\subsection*{Linear Polynomials Over $K$}

Towards understanding $m_L^{\sigma}(p(x)^t)$, let us first consider the simplest possible case, when $t = 1$ and $p(x) = x - \lambda$ is a linear polynomial with $\lambda \in K^{\times}$.

In this case we can say the following. Firstly, without going into detail, from the definition of $m_L^{\sigma}(x-\lambda)$ one has that for any $k \geq 1$, $m_L^{\sigma}(x-\lambda) \mid k$ if and only if there is $A \in \GL_k(L)$ with 
	\[
		A \sigma(A) \cdots \sigma^{r-1}(A) = \lambda I_k.
	\]
In particular, $m_L^{\sigma}(x-\lambda) = 1$ if and only if $\lambda$ is in the image of the norm map $N_{L/K} \colon L^{\times} \rightarrow K^{\times}$.

With a little more work (see Section \ref{sect:linearcase}), one can show that
\begin{itemize}
	\item $m_L^{\sigma}(x-\lambda)$ is the index of the cyclic central simple $K$-algebra $(L/K, \sigma, \lambda)$.
	\item If the period equals the index for elements of $\Br(K)$, then $m_L^{\sigma}(x-\lambda)$ is the order of $\lambda$ in 
	\[
		K^{\times} / N_{L/K}(L^{\times}).
	\]
\end{itemize}

We thus see that the semilinear Schur indices $m_L^{\sigma}(x-\lambda)$ for $\lambda \in K^{\times}$ are related to standard objects, and in many cases can be algorithmically determined, such as over local or global fields.

\subsection*{General Polynomials}

Our main result is a description of $m_L^{\sigma}(p(x)^t)$ for any $(p(x),t) \in \sP_L \times \bZ_{\geq 1}$. We do this by relating $m_L^{\sigma}(p(x)^t)$ to the linear case above, up to potentially enlarging the field.

More precisely, let $\alpha \in \overline{L}$ be any root of $p(x)$, let $K_{p(x)}$ denote the subfield of $L$ generated by $K$ and the coefficients of $p(x)$, and set $r(p) \coloneqq [K_{p(x)} : K]$.

In this setting, one can consider the field extension $L(\alpha) / K(\alpha)$, which is a cyclic Galois extension, with Galois group generated by $\tau_{\alpha}$, the automorphism of $L(\alpha)$ which acts trivially on $\alpha$ and by $\sigma^{r(p)}$ on $L$. Our main result can then be stated as follows.

\begin{theoremA}
	Suppose that $p(x) \in \sP_L$, $t \geq 1$, and $\alpha \in \overline{L}$ is a root of $p(x)$. Then
	\[
		m_L^{\sigma}(p(x)^t) = m_{L(\alpha)}^{\tau_{\alpha}}(x - \alpha).
	\]	
\end{theoremA}

In particular, we can use the results above for linear polynomials valued in $K$ to compute $m_L^{\sigma}(p(x)^t)$ for any $(p(x),t) \in \sP_L \times \bZ_{\geq 1}$. For example, if $L$ is a local or global field, Theorem A tells us that $m_L^{\sigma}(p(x)^t)$ is the order of $\alpha$ in 
\[
K(\alpha)^{\times} / N_{L(\alpha) / K(\alpha)}(L(\alpha)^{\times}),
\]
which can be determined using algorithms for solving norm equations \cite{AccKlu, Coh, LJMS, Denis}.

We prove Theorem A via the following series of steps.

\begin{enumerate}
	\item From general $p(x) \in \sP_L$ we reduce to $p(x)$ with coefficients in $K$. This follows a general principle that when computing semilinear Schur indices, one can always ``factor out'' the Galois action (Lemma \ref{lem:reductiontopolyoverK}).
	\item For $p(x) \in \sP_L$ with coefficients in $K$, we prove that $m_L^{\sigma}(p(x)) = m_{L(\alpha)}^{\tau_{\alpha}}(x - \alpha)$ by establishing a correspondence between semilinear extensions of the respective underlying linear representations over $L$ and $L(\alpha)$ (Theorem \ref{thm:schurindexpx}).
	\item Finally, we show that the irreducible case when $t = 1$ suffices in general, by a lifting argument which shows that $m_L^{\sigma}(p(x)^t) = m_L^{\sigma}(p(x))$ for any $t \geq 1$ (Theorem \ref{thm:mpequalsmpt}).
\end{enumerate}

\subsection*{Applications}

We obtain the following applications from the results described above.

\subsubsection*{A Normal Form for Matrices up to $\sigma$-Conjugacy}

For any $n \geq 1$ and $A, B \in \GL_n(L)$, we say that $A$ and $B$ are $\sigma$-conjugate, written $A \sim_{\sigma} B$, if there is some $P \in \GL_n(L)$ with $P^{-1} A \sigma(P) = B$.

There is a canonical identification between $\GL_n(L) / \sim_{\sigma}$ and isomorphism classes of $n$-dimensional objects of $\Rep_L^{\sigma}(\bZ)$. As a consequence, there are, for each orbit $\OO \in \Gal(L/K) \backslash (\sP_L \times \bZ_{\geq 1})$, invertible matrices $A_{\OO}$ over $L$ such that for any $A \in \GL_n(L)$, $A$ is $\sigma$-conjugate to the block sum
\[
A \sim_{\sigma} A_{\OO_1} \oplus \cdots \oplus A_{\OO_m}
\]
for some $\OO_1, ... , \OO_m \in \Gal(L/K) \backslash (\sP_L \times \bZ_{\geq 1})$, unique up to permutation.

In general, the representatives $A_{\OO}$ are unique only up to $\sigma$-conjugation, and depend on a choice. This is much like the theory of the Jordan normal form, or generalised Jordan normal form, where one chooses as representatives the blocks $J_t(\lambda)$ in the first case and companion matrices in the second. The difference now is that there is no such canonical choice in general. Furthermore, the fact that there is no formula for $A_{\OO}$ means that representatives can be hard to find. 

We do, however, give a recipe for constructing representatives $A_{\OO}$ in certain cases. This includes the case that $L/K$ has prime degree, or that $L/K$ is an unramified extension of local fields.

Given $A$, we also provide an algorithm to determine the $\OO_1, ... , \OO_m$, and so the block decomposition of $A$. This algorithm rests on the computation of some semilinear Schur indices.

As another consequence of properties of $\Rep_L^{\sigma}(\bZ)$, we show one can test $\sigma$-conjugacy in terms of ordinary conjugacy: for $A, B \in \GL_n(L)$, $A \sim_{\sigma} B$ if and only if $N_n^{\sigma}(A) \sim N_n^{\sigma}(B)$, where 
\[
N_n^{\sigma}(A) \coloneqq A \sigma(A) \cdots \sigma^{r-1}(A)
\]
In particular, one can test $\sigma$-conjugacy computationally using existing linear-algebraic algorithms.

\subsubsection*{A Description of the Category of Isocrystals Over a Finite Field}

For a perfect field $k$ of characteristic $p > 0$, we may consider the category of isocrystals over $k$, defined as follows.

Write $K(k) \coloneqq W(k)[1/p]$ for the field of fractions of the ring of Witt vectors of $k$, and $\varphi$ for the Frobenius automorphism of $K(k)$.

\begin{defn*}
The category $\Isoc(k)$ consists of pairs $(V,F)$, where $V$ is a finite-dimensional $K(k)$-vector space and $F \colon V \rightarrow V$ is a $\varphi$-linear bijection. Morphisms $(V_1, F_1) \rightarrow (V_2,F_2)$ are $K(k)$-linear maps $f \colon V_1 \rightarrow V_2$ such that $f \circ F_1 = F_2 \circ f$.
\end{defn*}

The category $\Isoc(k)$ plays an important role in arithmetic geometry. For example, taking the Dieudonn\'{e} module of a $p$-divisible group defines an equivalence from the category of $p$-divisible groups over $k$ up to isogeny to the full subcategory of $\Isoc(k)$ whose objects have all slopes in the interval $[0,1]$ (see \cite{Manin} when $k = \overline{k}$, or \cite{Marrama} for general $k$). Isocrystals also arise from cohomology: when $\bF_q$ is a finite field, then assuming the Tate conjecture, crystalline cohomology gives a fully faithful embedding from the category $\text{Mot}(\bF_q) \otimes \bQ_p$ of motives over $\bF_q$ with coefficients in $\bQ_p$ to $\Isoc(\bF_q)$ \cite[Prop.\ 1.17]{MilMot}.

When $k$ is algebraically closed, $\Isoc(k)$ has a well-known explicit description: from the Dieudonn\'{e}-Manin classification, $\Isoc(k)$ is a semisimple category, and the simple objects are explicitly defined and indexed by $\bQ$ \cite{Manin}. One other case when $\Isoc(k)$ is well understood is when $k = \bF_p$. In this case $K(\bF_p) = \bQ_p$, and $\varphi$ is the identity automorphism, hence $\Isoc(\bF_p)$ is simply the category of finite-dimensional $\bQ_p$-linear representations of $\bZ$. Unlike $\Isoc(\overline{\bF_p})$, $\Isoc(\bF_p)$ is not semisimple, and isomorphism classes of objects are in bijection with invertible matrices over $\bQ_p$ up to conjugacy.

In general, there is similarly a canonical identification
\[
\Isoc(k) \xrightarrow{\sim} \Rep_{K(k)}^{\varphi}(\bZ)
\]
and thus $\Isoc(k)$ may be interpreted as the category of semilinear representations of $\bZ$ over $K(k)$, where $\bZ$ acts on $K(k)$ via $\varphi$. In particular, $\Isoc(k)$ is amenable to the methods of this paper precisely when $\varphi$ has finite order, which is exactly when $k = \bF_q$ is a finite field. 

In this case, for $q = p^r$, $K(\bF_q) = \bQ_q$, the unique unramified extension of $\bQ_p$ of degree $r$, and $\varphi = \varphi_q$ is the lift of the Frobenius automorphism $x \mapsto x^p$ of $\Gal(\bF_q / \bF_p)$, of order $r$.

From our general analysis of semilinear representations of $\bZ$, we immediately obtain: 
\begin{itemize}
	\item Every object of $\Isoc(\bF_q)$ decomposes uniquely into indecomposable isocrystals.
	\item The set of indecomposable isocrystals over $\bF_q$ is in bijection with the set
	\[
		\left(\Gal_{\bQ_p} \!\! \backslash \overline{\bQ_p}^{\times} \right) \times \bZ_{\geq 1}
	\]
	where the subset of simple isocrystals corresponds to pairs with second factor equal to $1$.
	\item Any isocrystal is determined by its ``linearisation'': if $A_1, A_2 \in \GL_n(\bQ_q)$ are matrices defining the isocrystals $(\bQ_q^{n}, F_1)$ and $(\bQ_q^{n}, F_2)$, where $F_1(\mathbf{v}) \coloneqq A_1 \varphi(\mathbf{v})$ and $F_2(\mathbf{v}) \coloneqq A_2 \varphi(\mathbf{v})$, then $(\bQ_q^{n}, F_1) \cong (\bQ_q^{n}, F_2)$ if and only if $N^{\varphi_q}_n(A_1)$ and $N^{\varphi_q}_n(A_2)$ are conjugate in $\GL_n(\bQ_q)$. 
	\item An intrinsic description of the $B \in \GL_n(\bQ_q)$ arising as the linearisation of some isocrystal.
\end{itemize}

Here, each matrix $N^{\varphi_q}_n(A_i)$ records the action of the linear operator $F_i^r$, and the third point is the statement that this linear operator determines the isocrystal up to isomorphism.

For simple isocrystals the above bijection agrees with that of Kottwitz, who described the map 
\[
\Irr(\Isoc(\bF_q)) \xrightarrow{\sim} \Gal_{\bQ_p} \!\! \backslash \overline{\bQ_p}^{\times}
\]
which takes a simple isocrystal $(V,F)$ to a root of the minimal polynomial of $F^r$ \cite[Lem.\ 11.2]{Kott}.

As a consequence of our results, we describe the inverse of this map explicitly. The resulting construction of simple objects we give is significantly more involved than that of the Dieudonn\'{e}-Manin classification. We note that this is something that is not possible for $\Rep_L^{\sigma}(\bZ)$ over general base fields $L$, and rests on the fact that $\bQ_q / \bQ_p$ is an unramified extension of local fields.

To describe this inverse map, fix an algebraic closure $\overline{\bQ_p}$, and take $\alpha \in \overline{\bQ_p}^{\times}$.

\begin{enumerate}
\item Let $m_{\alpha, p}$, $m_{\alpha, q}$ denote the minimal polynomials of $\alpha$ over $\bQ_p$ and $\bQ_q$ respectively, and set
\[
r_{\alpha} \coloneqq \frac{\deg(m_{\alpha, p})}{\deg(m_{\alpha, q})}.
\]
\item Write $\alpha = \pi^s u$, where $\pi$ is a uniformiser of $\bQ_{p}(\alpha)$ and $u \in \OO_{\bQ_p(\alpha)}^{\times}$.
\item Let $v \in \OO_{\bQ_q(\alpha)}^{\times}$ with $N_{\bQ_q(\alpha) / \bQ_p(\alpha)}(v) = u$, possible because $\bQ_q(\alpha) / \bQ_p(\alpha)$ is unramified. 
\item Let $s_{\alpha} \coloneqq \frac{r}{r_{\alpha}} \in \bZ_{\geq 1}$, $g_{\alpha} \coloneqq \gcd(s_{\alpha},s)$, and $m_{\alpha} \coloneqq \frac{s_{\alpha}}{g_{\alpha}}$, and set
\[
B_{\alpha} \coloneqq v \cdot \begin{pmatrix}
    0 & 0 &\cdots & 0 & \pi^{\frac{s}{g_{\alpha}}} \\
    1 & 0 & \cdots & 0 & 0 \\
    \vdots & \vdots & \ddots & \ddots & \vdots \\
    0 & 0 & \cdots & 0 & 0 \\
    0 & 0 & \cdots & 1 & 0
\end{pmatrix} \in \GL_{m_{\alpha}}(\bQ_q(\alpha)).
\]
\item Let $n_{\alpha} \coloneqq \deg(m_{\alpha, q})$, and write
\[
M_{\alpha} \colon M_{m_{\alpha}}(\bQ_q(\alpha)) \hookrightarrow M_{m_{\alpha} n_{\alpha}}(\bQ_q)
\]
for the map which on each entry $\beta \in \bQ_q(\alpha)$ is the natural inclusion 
\[
\bQ_q(\alpha) \hookrightarrow M_{n_{\alpha}}(\bQ_q)
\]
sending $\beta$ to the matrix of left multiplication by $\beta$ with respect to the basis $1, \alpha, \cdots , \alpha^{n_{\alpha}-1}$.
\item Set $d_{\alpha} \coloneqq m_{\alpha} n_{\alpha} r_{\alpha}$ and $A_{\alpha}$ to be the block matrix
\[
A_{\alpha} \coloneqq \begin{pmatrix}
    0 & 0 &\cdots & 0 & M_{\alpha}(B_{\alpha}) \\
    I_{m_{\alpha} n_{\alpha}} & 0 & \cdots & 0 & 0 \\
    \vdots & \vdots & \ddots & \ddots & \vdots \\
    0 & 0 & \cdots & 0 & 0 \\
    0 & 0 & \cdots & I_{m_{\alpha} n_{\alpha}} & 0
\end{pmatrix} \in \GL_{d_{\alpha}}(\bQ_q).
\]
\item Finally, let $W_{\alpha}$ denote the isocrystal $W_{\alpha} \coloneqq (\bQ_q^{d_{\alpha}}, F_{\alpha})$, where
\[
F_{\alpha}(\mathbf{v}) \coloneqq A_{\alpha} \varphi(\mathbf{v}).
\]
\end{enumerate}

Then $W_{\alpha}$ is simple, and (up to isomorphism) is independent of the choices involved above. Furthermore, $W_{\alpha}$ only depends on the $\Gal_{\bQ_p}$-orbit of $\alpha$, and this induces the identification
\[
\Gal_{\bQ_p} \!\! \backslash \overline{\bQ_p}^{\times} \xrightarrow{\sim} \Irr(\Isoc(\bF_q)).
\]
We similarly provide an explicit construction of all indecomposable isocrystals, see Section \ref{sect:isocrystals}.

For example, this shows that the isocrystal of an abelian variety over $\bF_q$, which determines its $p$-divisible group up to isogeny, can be written down explicitly from its Weil polynomial.

As one application of our description of $\Isoc(\bF_q)$, we study how isocrystals behave under base change. For example, we classify which isocrystals over $\bF_q$ become isomorphic to their Dieudonn\'{e}-Manin description already over some finite extension $\bF_{q'} / \bF_q$ (meaning in some $\Isoc(\bF_{q'})$).

\subsection*{Further Directions} Our work raises many further questions that merit further study.

\subsubsection*{Higher Rank Norm Equations}

For the category of isocrystals, $\Rep_{\bQ_q}^{\varphi_q}(\bZ)$, we are able to give explicit models for all indecomposable representations. In general, for $\Rep_{L}^{\sigma}(\bZ)$ this is equivalent to the question of explicitly constructing the representatives $A_{\OO}$ for the semilinear normal form, and requires solving equations of the type
\[
A \sigma(A) \cdots \sigma^{r-1}(A) = \lambda I_n
\]
for $A \in \GL_n(L)$, where $n \geq 1$ and $\lambda \in K^{\times}$. These equations are solvable when $n = 1$ and $L$ is a number field, as these equations are norm equations for which there exists a developed computational theory \cite{Coh, Denis}. It is thus natural to ask if over number fields these algorithms can be generalised to higher rank, as this would give an explicit description of $\Rep_{L}^{\sigma}(\bZ)$ in this case.

\subsubsection*{Non-Invertible Matrices}

Here we have classified matrices in $\GL_n(L)$ up to $\sigma$-conjugacy, but it is natural to ask how this extends to $M_n(L)$. In the special case that $L$ is $\bC$, and $\sigma$ is complex conjugation, this was studied in \cite{HongHorn}. For example the fact that $A \overline{A} \sim B \overline{B}$ forces $A \sim_{\sigma} B$ does not hold for non-invertible matrices in general, and they give an extra condition on the rank which allows one to extend this property from $\GL_n(\bC)$ to all of $M_n(\bC)$.

\subsubsection*{Integral Structures}

Here we have classified isocrystals over $\bF_q$, and it is natural to ask for an explicit classification of the $W(\bF_q)$-lattices stable under $F$ and $p F^{-1}$ inside a given isocrystal up to isomorphism. For slopes in $[0,1]$, this amounts to classifying $p$-divisible groups over $\bF_q$ up to isomorphism rather than isogeny. For Dieudonn\'{e} modules of abelian varieties, explicit classifications are known in many cases \cite{BKM, WaterhouseAV}, but not in general. In the proof of \cite[Thm.\ 5.1]{WaterhouseAV}, such lattices are classified for principal abelian varieties (those which are simple with endomorphism ring the maximal order of a field) using an explicit model of the isocrystal similar to ours in the case $m_{\alpha} = 1$, and so we expect the explicit models given here to be useful in the case that $m_{\alpha} > 1$.

\subsection*{Organisation of the Paper}

In Section \ref{sect:GJNF} we recall the relationship between matrices, $L[x]$-modules, and representations of $\bZ$, and the generalised Jordan normal form. In Section \ref{sect:semilinearbackground} we explain how this translates to the semilinear setting. In Section \ref{sect:schurindicesresults} we prove Theorem A. In Section \ref{sect:representatives} we consider the question of constructing explicit representatives for the semilinear normal form, and in Section \ref{sect:basechange} describe how semilinear representations of $\bZ$ behave under base change. In Section \ref{sect:localfields} we specialise to the case of local fields, and prove our results about the category of isocrystals above. In Section \ref{sect:cyclicgroups}, we give some further results for semilinear representations of finite cyclic groups, and finally in Section \ref{sect:algorithms} we summarise the algorithms appearing in the text, and include implementations over number fields.

\subsection*{Notation}

Throughout, $L$ denotes a field, and $\sigma$ an automorphism of $L$ of finite order $r \geq 1$. We write $K$ for the fixed field of $\sigma$, so $L/K$ is Galois with $\Gal(L/K) = \langle \sigma \rangle$, cyclic of order $r$.

We write $\sP_L$ for the set of monic irreducible polynomials over $L$ with non-zero constant term, and for $p(x) \in \sP_L$, write $K_{p(x)}$ for the subfield of $L$ generated by $K$ and the coefficients of $p(x)$ and denote $r(p) \coloneqq [K_{p(x)}:K]$. For $\alpha \in \overline{L}$ a root of $p(x)$, $\tau_{\alpha} \in \Gal(L(\alpha) / K(\alpha))$ denotes the automorphism of $L(\alpha)$ which acts by $\sigma^{r(p)}$ on $L$ and fixes $\alpha$.

\subsection*{AI Usage Statement} The main ideas and results of this paper were all established without the use of LLMs. Free versions of LLMs were used to check for errors, write the associated code, and suggested simplifications for some of the proofs. The paper was written without the use of LLMs.

\subsection*{Acknowledgements} Both authors would like to thank Rhiannon Savage and Niki Kalaydzhieva for organising the AFFIRM research school for women and other underrepresentated genders in mathematics, where this project originated. During this work, the first author was supported by the departmental grant Progetto Sviluppo Dipartimentale - UNIPD PSDIP23O88 at the University of Padova, and the Heilbronn Institute of Mathematical Research at the University of Bristol.

\section{Recollections on the Generalised Jordan Normal Form}\label{sect:GJNF}

Let $L$ be a field. In this section we recall the relationship between $L[x]$-modules, matrices over $L$ up to conjugation, and $L$-linear representations of $\bZ$.

\begin{itemize}
	\item Given any matrix $A \in M_n(L)$, there is an associated $L[x]$-module $L^n_A$, which is $L^n$ as an $L$-vector space, and action of $x$ via left multiplication by $A$.
	\item Every $L[x]$-module which is finite-dimensional over $L$ is isomorphic to some $L^n_A$.
	\item Given $A, B \in M_n(L)$, we have that $L^n_A \cong L^n_B$ as $L[x]$-modules if and only if there is some $P \in \GL_n(L)$ with $A P = P B$.
	\item This induces an identification, for any $n \geq 1$, from the set of $n \times n$-matrices over $L$ up to conjugation, to the set of $n$-dimensional $L[x]$-modules up to isomorphism.
	\item Given $f(x) \in L[x]$, of degree $n \geq 1$, one can consider the cyclic $L[x]$-module
	\[
	M_f \coloneqq L[x] / (f(x)).
	\]
	\item The companion matrix of $f(x) = a_n x^n + \cdots a_1 x + a_0$, is the matrix 
\[
C(f) = \begin{pmatrix}
	0 & 0 &\cdots & 0 & -a_0 / a_n \\
	1 & 0 & \cdots & 0 & -a_1 / a_n \\
	\vdots & \vdots & \ddots & \ddots & \vdots \\
	0 & 0 & \cdots & 0& -a_{n-2} / a_n  \\
	0 & 0 & \cdots & 1 & -a_{n-1} / a_n
\end{pmatrix},
\]
which records the action of $x$ with respect to the $L$-basis $1, x, ... , x^{n-1}$ of $M_f$.
\item In particular, there is a canonical isomorphism
\[
M_f \cong L^n_{C(f)},
\]
where the $L$-basis $1, x, ... , x^{n-1}$ of $M_f$ of $M_f$ corresponds to the the standard basis of $L^n_{C(f)}$.
\end{itemize}

From the classification of finitely generated modules over PIDs, we have the following.

\begin{itemize}
\item Any $L[x]$-module which is finite-dimensional over $L$ is isomorphic to
\[
M_{p_1^{r_1}} \oplus \cdots \oplus M_{p_m^{r_m}}
\]
for $p_1(x), ... , p_m(x)$ irreducible monic polynomials\footnote{which need not be distinct} over $L$, and $r_1, ... , r_m \geq 1$.

\item This description is unique up to reordering of the factors.

\item For $p(x) \in L[x]$ monic and irreducible, and $t \geq 1$, the $L[x]$-module $M_{p^t}$ is indecomposable, and irreducible if and only if $t = 1$.
\end{itemize}

We may translate this to a statement about matrices to obtain a classification of matrices up to conjugation over $L$, known as the generalised Jordan Normal Form, or GJNF for short. Here the symbol $\oplus$ denotes the block sum of matrices.

\begin{itemize}
\item Any $A \in M_n(L)$ is conjugate to a matrix of the form
\[
C(p_1^{r_1}) \oplus \cdots \oplus C(p_m^{r_m})
\]
for $p_1(x), ... , p_m(x)$ irreducible monic polynomials over $L$, and $r_1, ... , r_m \geq 1$.

\item This description is unique up to reordering of the blocks.
\end{itemize}

We note that this is indeed a generalisation of the classical Jordan norm (which applies when $L$ is algebraically closed), as for $\lambda \in L$ and $t \geq 1$, the matrices $J_t(\lambda)$ and $C((x-\lambda)^t)$ are conjugate.

In practice, one can compute the GJNF as follows of some $A \in M_n(L)$ as follows.

\begin{enumerate}
	\item Factor the characteristic polynomial $\chi_A(x) = q_1(x)^{s_1} \cdots q_l(x)^{s_l}$.
	\item Each $p_i(x)$ will be one of these $q_1(x), ... , q_l(x)$.
	\item For $1 \leq i \leq l$, set $d_0(q_i) \coloneqq 0$, and for $k \geq 1$,
	\[
	d_k(q_i) \coloneqq \dim_L \ker(q_i^k(A)).
	\]
	\item The number times $q_i$ appears in the decomposition of $A$ with exponent $\geq k$ is
\[
\frac{d_k(q_i) - d_{k-1}(q_i)}{\deg(q_i(x))}.
\]
\item Therefore, the number of times the block $C(p_i(x)^k)$ appears in the GJNF of $A$ is
\[
\frac{2 d_k(q_i) - d_{k+1}(q_i) - d_{k-1}(q_i)}{\deg(q_i(x))}.
\]
\end{enumerate}

How do $L[x]$-modules and matrices relate to representations of $\bZ$? Let $\Rep_L(\bZ)$ denote the category of finite-dimensional $L$-linear representations of $\bZ$.

\begin{itemize}
	\item The group ring $L[\bZ]$ of $\bZ$ over $L$ is canonically identified with $L[x,x^{-1}]$, where $x$ corresponds to $1 \in \bZ$.
	
	\item Thus $\Rep_L(\bZ)$ is canonically identified with the category of $L[x,x^{-1}]$-modules which are finite-dimensional over $L$, which itself is canonically identified with the full subcategory of $L[x]$-modules which are finite-dimensional over $L$ for which $x$ acts invertibly.
	
	\item From a matrix $A \in \GL_n(L)$ one obtains a $\bZ$-representation on $L^n$ where $1 \in \bZ$ acts by $A$, and the corresponding $L[x,x^{-1}]$-module is $L^n_A$.
	
	\item This induces an identification, for any $n \geq 1$ from $\GL_n(L)$ up to conjugation to the set of $n$-dimensional objects of $\Rep_L(\bZ)$ up to isomorphism.
\end{itemize}

In terms of the GJNF, $A$ is invertible (and so corresponds to a $\bZ$-representation) if and only if none of the $p_i(x)$ appearing in the block decomposition of $A$ are equal to $x$. 

\begin{defn}
	$\sP_{L}$ is the set of monic irreducible polynomials $p(x) \in L[x]$ with $p(x) \neq x$.
\end{defn}

Let $\Irr_L(\bZ)$ and $\Ind_L(\bZ)$ denote the sets of isomorphism classes of irreducible and indecomposable finite-dimensional representations of $\bZ$ over $L$. For $p(x) \in \sP_L$ and $t \geq 1$, let $V_{p(x),t}$ denote the $\bZ$-representation where the generator $1 \in \bZ$ acts via the matrix $C(p^t)$. As a direct consequence of the above discussion, we have the following, which we record for later use.

\begin{prop}\label{prop:indZreps}
	The map
	\[
	\sP_L \times \bZ_{\geq 1} \rightarrow \Ind_L(\bZ), \qquad (p(x),t) \mapsto V_{p(x),t}
	\]
	is a bijection. The subset $\sP_L \equiv \sP_L \times \{1\} \subset \sP_L \times \bZ_{\geq 1}$ corresponds to the subset $\Irr_L(\bZ) \subset \Ind_L(\bZ)$.
\end{prop}

Finally, we also record here the following basic lemma, which we will use multiple times.

For $M \in M_n(L)$, let $L[M] \subset M_n(L)$ denote the $L$-subalgebra generated by $M$, and let $\fC(M)$ denote the centraliser of $M$ in $M_n(L)$.

\begin{lemma}\label{lem:centralisercyclicmodule}
	Suppose that $f(x) \in L[x]$ has degree $n \geq 1$. Then $\fC(C(f)) = L[C(f)]$, and the natural map
	\[
	L[x] / f(x) \rightarrow L[C(f)], \qquad x \mapsto C(f)
	\]
	is an $L$-algebra isomorphism.
\end{lemma}

\begin{proof}
	The given map is clearly an isomorphism, as $C(f)$ has minimal polynomial $f(x)$. Furthermore, $L[C(f)] \subset \fC(C(f))$, and so it suffices to show the reverse inclusion.
	
	Suppose then that $B \in \fC(C(f))$, and let $e_1, ... , e_n$ denote the standard basis of $L^n$. For any $0 \leq i \leq n-1$,
	\[
	C(f)^i e_1 = e_{i+1},
	\]
	hence $\sB \coloneqq \{C(f)^i e_1\mid 0 \leq i \leq n - 1\}$ forms an $L$-basis of $L^n$, and we may write
	\[
	B e_1 = \sum_{i = 0}^{nt - 1} a_i C(f)^i e_1
	\]
	for some $a_i \in L$. We have that for any $0 \leq j \leq n - 1$,
	\[
		B (C(f)^j e_1) = C(f)^j B e_1 =  C(f)^j \left(\sum_{i = 0}^{n - 1} a_i C(f)^{i} e_1 \right) = \left(\sum_{i = 0}^{n - 1} a_i C(f)^i\right) C(f)^j e_1,
	\]
	and therefore
	\[
	B = \sum_{i = 0}^{n - 1} a_i C(f)^i
	\]
	because $\sB$ form an $L$-basis of $L^{n}$.
\end{proof}

\section{Semilinear Representations of $\bZ$ and a Semilinear Normal Form for Matrices}\label{sect:semilinearbackground}

In this paper we are interested in semilinear representations of $\bZ$, arising as follows. Suppose henceforth that $\sigma$ is an automorphism of $L$ of finite order $r \geq 1$, and set $K \coloneqq L^{\sigma}$. The extension $L/K$ is Galois, with Galois group $\langle \sigma \rangle$, and we let $\bZ$ act on $L$ via the homomorphism
\[
\tilde{\sigma} \colon \bZ \rightarrow \Gal(L/K), \qquad \tilde{\sigma}(m) = \sigma^m.
\]

In this context, we can consider the category of semilinear representations of $\bZ$ over $L$, which we denote by $\Rep_L^{\sigma}(\bZ)$. As $\bZ$ is cyclic, one can think of this category as follows.

Objects of $\Rep_L^{\sigma}(\bZ)$ are pairs $(V,T)$, where $V$ is a finite-dimensional vector space over $L$, and $T \colon V \rightarrow V$ is an invertible $\sigma$-linear endomorphism of $V$, and a morphism $f \colon (V_1, T_1) \rightarrow (V_2, T_2)$ is an $L$-linear map $f \colon V_1 \rightarrow V_2$ such that $f \circ T_1 = T_2 \circ f$. 

\subsection{Module Theoretic and Matrix Descriptions}
As in Section \ref{sect:GJNF}, when the action of $\bZ$ on $L$ is trivial, there are two other points of view on $\Rep_L^{\sigma}(\bZ)$ in terms of both modules and matrices. 

\begin{defn}
	Let $L_{\sigma}[x,x^{-1}]$ be the ring which is $L[x,x^{-1}]$ as an abelian group, but with multiplication determined by 
	\[
	\lambda x^m \cdot \mu x^n \coloneqq \lambda \sigma^m(\mu) x^{m+n}
	\]
	for $\lambda, \mu \in L$ and $m,n \in \bZ$.
\end{defn}

There is an equivalence of categories between $\Rep_L^{\sigma}(\bZ)$ and the category of $L_{\sigma}[x,x^{-1}]$-modules which are finite-dimensional over $L$, which sends a pair $(V,T)$ to the $L$-vector space $V$ with action of $x \cdot v \coloneqq T(v)$, and is the identity on morphisms.

Given $A \in \GL_n(L)$, one can form the $L_{\sigma}[x,x^{-1}]$-module $L^{n ,\sigma}_A$, which is $L^n$ as an $L$-vector space, with action of $x$ via the map
\[
x \cdot \mathbf{v} \coloneqq A \sigma(\mathbf{v}).
\]
For example, for $m \geq 1$ with have that
\begin{equation}\label{eqn:higherpowerofxaction}
x^m \cdot \mathbf{v} = A \sigma(A) \cdots \sigma^{m-1}(A) \sigma^m(\mathbf{v}), \qquad \text{and} \qquad x^{-1} \cdot \mathbf{v} = \sigma^{-1}(A^{-1} \mathbf{v}).
\end{equation}

Any $L_{\sigma}[x,x^{-1}]$-module that is finite-dimensional over $L$ is of the form $L^{n ,\sigma}_A$ for some $A \in \GL_n(L)$, and given $A, B \in \GL_n(L)$, we have that $L^{n, \sigma}_A \cong L^{n, \sigma}_B$ as $L_{\sigma}[x,x^{-1}]$-modules if and only if there is some $P \in \GL_n(L)$ with $A \sigma(P) = P B$. This motivates the following definition.

\begin{defn}
	For $A, B \in \GL_n(L)$, we say that $A$ and $B$ are $\sigma$-conjugate, written $A \sim_{\sigma} B$, if there exists $P \in \GL_n(L)$ with $P^{-1} A \sigma(P) = B$.
\end{defn}

In particular for any $n \geq 1$ the map $A \mapsto L^{n, \sigma}_A$ induces an identification from $\GL_n(L) / \!\! \sim_{\sigma}$ to the set of $n$-dimensional $L_{\sigma}[x,x^{-1}]$-modules up to isomorphism, or equivalently $n$-dimensional objects of $\Rep_L^{\sigma}(\bZ)$ up to isomorphism.

\subsection{The Restriction Functor and Semilinear Schur Indices}\label{sect:backgroundsemilinearreps}

We are interested in describing $\Rep_L^{\sigma}(\bZ)$, so let us now describe how the results of \cite{TAYSchur} apply in this setting. For the subgroup $r \bZ \leq \bZ$, we have a natural restriction functor
\[
(-)|_{r \bZ} \colon \Rep_L^{\sigma}(\bZ) \rightarrow \Rep_L(r \bZ),
\]
noting that $r\bZ$ acts on $L$ trivially because $\sigma$ has order $r$. The first main result of \cite{TAYSchur} is as follows. 

\begin{thm}[{\cite[Theorem A]{TAYSchur}}]\label{thm:isomonrestrictions}
	If $V_1, V_2 \in \Rep_L^{\sigma}(\bZ)$ with $V_1|_{r\bZ} \cong V_2|_{r \bZ}$, then $V_1 \cong V_2$.
\end{thm}

In other words, Theorem \ref{thm:isomonrestrictions} states that there is at most one extension of a linear representation of $r \bZ$ to a semilinear representation of $\bZ$. Therefore, one can obtain a description of the objects of $\Rep_L^{\sigma}(\bZ)$ by describing exactly which objects of $\Rep_L(r \bZ)$ extend to $\bZ$, or equivalently by describing the essential image of the functor $(-)|_{r \bZ}$. How can we determine this image?

Any object of $\Rep_L^{\sigma}(\bZ)$ decomposes uniquely into a direct sum of indecomposable objects \cite[Lem.\ 2.14]{TAYSchur}, and one is reduced to understand the restrictions of $\Ind_L^{\sigma}(\bZ)$ and $\Irr_L^{\sigma}(\bZ)$, the sets of isomorphism classes of indecomposable and irreducible objects of $\Rep_L^{\sigma}(\bZ)$.

There is a natural left-action of $\Gamma \coloneqq \Gal(L/K)$ on the sets $\Ind_L(r \bZ)$ and $\Irr_L(r \bZ)$, for which $\gamma*(V,T)$ is the pair $({}^{\gamma} V, T)$, where ${}^{\gamma} V = V$ but with $L$-vector-space structure $\lambda * v \coloneqq \gamma^{-1}(\lambda)v$. We write $\Gamma_V$ for the stabiliser of $V$ in $\Gamma$.

\begin{thm}[{\cite[Theorem B]{TAYSchur}}]\label{thm:descofrestriction}
	There is a bijection
	\[
		\Ind_L^{\sigma}(\bZ) \xrightarrow{\sim} \Gamma \: \backslash \Ind_{L}(r \bZ)
	\]
	where $W \in \Ind_L^{\sigma}(\bZ)$ maps to any indecomposable direct summand $V \in \Ind_L(r \bZ)$ of $W|_{r \bZ}$.
	
	Furthermore, for any such $V$, there is some $m_L^{\sigma}(V) \in \bZ_{\geq 1}$ with,
	\[
		W|_{r \bZ} \cong \left( \bigoplus_{\gamma \in \Gamma / \Gamma_V} \gamma *V \right)^{m_L^{\sigma}(V)}.
	\]
\end{thm}

The number $m_K^L(V)$ is called the semilinear Schur index of $V$. These depend only on the $\Gamma$-orbit of $V$, and are hard to compute in general. A knowledge of the semilinear Schur indices $m_K^L(V)$ for all $V \in \Ind_L(r \bZ)$ allows one to give a description of the category $\Rep_L^{\sigma}(\bZ)$ in terms of $\Ind_{L}(r \bZ)$ and its action of $\Gamma$: 
\begin{itemize}
\item For any $V \in \Ind_L(r \bZ)$, $m_L^{\sigma}(V)$ copies of the orbit sum
\[
\bigoplus_{\gamma \in \Gamma / \Gamma_V} \gamma * V
\]
extends uniquely to a semilinear representation $\widetilde{V}$ of $\bZ$,
\item For $V_1, V_2 \in \Ind_L(r \bZ)$, $\widetilde{V_1} \cong \widetilde{V_2}$ if and only if $V_1$ and $V_2$ are in the same $\Gamma$-orbit,
\item All elements of $\Ind_{L}^{\sigma}(\bZ)$ arise this way.
\end{itemize}

We thus see that to describe $\Rep_L^{\sigma}(\bZ)$, we need to understand the $\Gamma$-set $\Ind_L(r \bZ)$, and the numbers $m_L^{\sigma}(V)$.

Let us first consider $\Ind_L(r \bZ)$. From Proposition \ref{prop:indZreps}, there is a bijection
\[
\sP_L \times \bZ_{\geq 1} \rightarrow \Ind_L(r \bZ), \qquad (p(x),t) \mapsto V_{p(x),t,r},
\]
where we write $V_{p(x),t,r}$ for the $r \bZ$-representation $L^{\deg(p)t}$, where the generator $r$ of $r \bZ$ acts on $L^{\deg(p)t}$ via the matrix $C(p^t)$. The representation $V_{p(x), t, r}$ is irreducible if and only if $t = 1$. 

This bijection is $\Gamma$-equivariant with respect to the natural action of $\Gal(L/K)$ on $\sP_L$ and the trivial action on $\bZ_{\geq 1}$, and we thus obtain a canonical identification
\[
\Gamma \: \backslash (\sP_L \times \bZ_{\geq 1}) \leftrightarrow \Ind_L^{\sigma}(\bZ).
\]

\begin{defn}
    For $p(x) \in \sP_L$ and $t \geq 1$, we write $W_{p(x),t} \in \Ind_{L}^{\hspace{0.05em} \sigma}(\bZ)$ for the indecomposable representation corresponding to $V_{p(x),t,r} \in \Ind_L(r\bZ)$. For an orbit $\OO \in \Gamma \: \backslash (\sP_L \times \bZ_{\geq 1})$, we set $W_{\OO} \coloneqq W_{p(x),t}$ for any representative $(p(x),t) \in \sP_L \times \bZ_{\geq 1}$.
\end{defn}

\begin{defn}\label{def:schurindexnotation}
	For $p(x) \in \sP_L$ and $t \geq 1$, we denote
	\[
	m_L^{\sigma}(p(x)^t) \coloneqq m_L^{\sigma}(V_{p(x), t, r}) \in \bZ_{\geq 1},
	\]
	the semilinear Schur index of $V_{p(x), t, r} \in \Ind_L(r\bZ)$ with respect to the action of $\bZ$ on $L$ via $\sigma$.
\end{defn}

To describe $\Rep_L^{\sigma}(\bZ)$, we are thus reduced to determining the numbers $m_L^{\sigma}(p(x)^t)$ for $p(x) \in \sP_L$. This is the topic of Section \ref{sect:schurindicesresults}.

\subsection{Matrix Interpretation}\label{sect:matrixsemilinearpov}

Before we consider the $m_L^{\sigma}(p(x)^t)$, let us explain how description of $\Rep_L^{\sigma}(\bZ)$ translates into statements about matrices, and a semilinear normal form.

\begin{defn}
    For $m \geq 1$, we write
    \[
    N_m^{\sigma} \colon M_m(L) \rightarrow M_m(L), \qquad N_{m}^{\sigma}(A) \coloneqq A \sigma(A) \cdots \sigma^{r-1}(A).
    \]
\end{defn}

\begin{remark}
    When $m = 1$, then by definition $N_{1}^{\sigma}(a) = N_{L/K}(a)$ for $a \in L^{\times}$, where
    \[
    N_{L/K} \colon L^{\times} \rightarrow K^{\times}
    \]
    denotes the norm map. Many of the nice properties of the norm $N_{L/K}$ do not carry over, however: $N_m^{\sigma}(M_m(L))$ need not contained in $M_m(K)$, $N_{m}^{\sigma}$ does not induce a group homomorphism $\GL_m(L) \rightarrow \GL_m(L)$, and $N_{m}^{\sigma}(A)$ depends on the order of the terms in the product.
\end{remark}

The map $N_m^{\sigma}(-)$ is of interest for the following reason. As in Section \ref{sect:GJNF}, for $A \in \GL_n(L)$ let $L^{n,r}_{A} \in \Rep_L(r \bZ)$ denote $L^n$ with action of generator $r \in r \bZ$ by left-multiplication by $A$.

\begin{lemma}\label{lem:explictrestriction}
	For $A \in \GL_n(L)$, $L^{n ,\sigma}_A|_{r \bZ} = L^{n,r}_{N_n^{\sigma}(A)}$.
\end{lemma}

\begin{proof}
	From the formula \eqref{eqn:higherpowerofxaction}, we have that the generator $r \in r \bZ$, which corresponds to $x^r \in L_{\sigma}[x,x^{-1}]$, acts by
	\[
	r \cdot \mathbf{v} = A \sigma(A) \cdots \sigma^{r-1}(A) \sigma^r(\mathbf{v}) = N_n^{\sigma}(A) \mathbf{v}
	\]
	because $\sigma^r = 1$, and thus $L^{n ,\sigma}_A|_{r \bZ}$ literally equals $L^{n,r}_{N_n^{\sigma}(A)}$.
\end{proof}

The translation of Theorem \ref{thm:isomonrestrictions} is as follows, which states that we can detect semilinear conjugacy after applying $N_n^{\sigma}(-)$.

\begin{cor}\label{cor:linsimimpliessim}
	If $A, B \in \GL_n(L)$ for some $n \geq 1$, with $N_n^{\sigma}(A) \sim N_n^{\sigma}(B)$, then $A \sim_{\sigma} B$.
\end{cor}

\begin{proof}
	We have that $N_n^{\sigma}(A) \sim N_n^{\sigma}(B)$ if and only if $L^{n,r}_{N_n^{\sigma}(A)} \cong L^{n,r}_{N_n^{\sigma}(B)}$ in $\Rep_L(r \bZ)$, or equivalently that $L^{n ,\sigma}_A|_{r \bZ} \cong L^{n ,\sigma}_B|_{r \bZ}$ in $\Rep_L(r \bZ)$ by Lemma \ref{lem:explictrestriction}. By Theorem \ref{thm:isomonrestrictions}, this implies that $L^{n ,\sigma}_A \cong L^{n ,\sigma}_B$ in $\Rep_L^{\sigma}(\bZ)$, and so $A \sim_{\sigma} B$.
\end{proof}

We can also immediately give a normal form for $\GL_n(L)$ up to $\sigma$-conjugacy as follows.

Note that in the generalised Jordan normal form of Section \ref{sect:GJNF}, a fixed matrix description is chosen for each indecomposable representation of $\bZ$: to the indecomposable module $V_{p(x),t}$, one chooses the matrix $C(p^t)$. Similarly, in the Jordan normal form over algebraically closed fields, one instead chooses $J_t(\lambda)$ as the matrix representative of the indecomposable module $V_{(x - \lambda), t}$.

We need to make a similar choice here. However, one cannot give an explicit formula for the matrix purely in terms of $p(x)$ and $t$. Instead, we work with any fixed set of representatives. We will return to the question of how to find such matrices and and their canonicity in Section \ref{sect:representatives}. 

From above, the indecomposable modules are $W_{p(x),t}$, indexed by $(p(x),t) \in \Gamma \: \backslash (\sP_L \times \bZ_{\geq 1})$.

\begin{defn}
	For each $\OO \in \Gamma \: \backslash (\sP_L \times \bZ_{\geq 1})$, fix some $A_{\OO} \in \GL_n(L)$ with
	\[
	W_{\OO} \cong L^{n, \sigma}_{A_{\OO}}.
	\]
\end{defn}

Then from the unique decomposition into indecomposables in $\Rep_L^{\sigma}(\bZ)$, we have the following.

\begin{cor}\label{cor:SNF}
For any $n \geq 1$, any $A \in \GL_n(L)$ is $\sigma$-conjugate to a block matrix
\[
A \sim_{\sigma} A_{\OO_1} \oplus \cdots \oplus A_{\OO_m}
\]
for $\OO_1, ... , \OO_m \in \Gamma \: \backslash (\sP_L \times \bZ_{\geq 1})$, where the symbol $\oplus$ denotes the block sum of matrices.

This description is unique up to reordering of the blocks.
\end{cor}

Now, given $A \in \GL_n(L)$, how can we compute the $\OO_1, ... , \OO_m \in \Gamma \: \backslash (\sP_L \times \bZ_{\geq 1})$ appearing in its semilinear normal form? 

\begin{enumerate}
	\item 
First, one computes $N^{\sigma}_n(A)$, which corresponds to restriction to the subgroup $r \bZ$, and finds the GJNF of $N^{\sigma}_n(A)$ using the algorithm described in Section \ref{sect:GJNF}, so
\[
N^{\sigma}_n(A) = C(q_1^{r_1}) \oplus \cdots \oplus C(q_m^{r_m})
\]
for some $q_1, \cdots, q_m \in \sP_L$ and $r_1, ... , r_m \geq 1$.
\item Because $N^{\sigma}_n(A)$ corresponds to the restriction of a semilinear representation of $\bZ$, we have that for any $\gamma \in \Gamma$ the number of copies of any $C(q^r)$ occurring equals the number of copies of $C(\gamma * q^r)$.

\item We may thus relabel and write 
\[
N^{\sigma}_n(A) = \left( \bigoplus_{(q,t) \in \OO_1} C(q^{t}) \right)^{\oplus k_1} \oplus \cdots \oplus \left( \bigoplus_{(q,t) \in \OO_k} C(q^{t})\right)^{\oplus k_m}
\]
for distinct orbits $\OO_1, ... , \OO_m \in \Gamma \: \backslash (\sP_L \times \bZ_{\geq 1})$.

\item The semilinear normal form of $A$ is then
\[
A \sim_{\sigma} A_{\OO_1}^{\oplus \frac{k_1}{m^{\sigma}_L(q_1(x)^{r_1})}} \oplus \cdots \oplus A_{\OO_m}^{\oplus \frac{k_m}{m^{\sigma}_L(q_m(x)^{r_m})}},
\]
where for $1 \leq i \leq m$, $(q_i,r_i)$ is any representative for the orbit $\OO_i$.
\end{enumerate}

In particular, we see that a knowledge of the semilinear Schur indices $m^{\sigma}_L(q_i(x)^{r_i})$ is needed to algorithmically find the semilinear normal form of elements of $\GL_n(L)$.

\section{Semilinear Schur Indices for Representations of $\bZ$}\label{sect:schurindicesresults}

In this section we study semilinear Schur indices for the group $\bZ$. Our strategy to describe $m_L^{\sigma}(p(x)^t)$ for general $p(x) \in \sP_L$ and $t \geq 1$ is via the following series of steps.
\begin{enumerate}
	\item From general $p(x) \in \sP_L$ we reduce to $p(x)$ with coefficients in $K$ (Lemma \ref{lem:reductiontopolyoverK}).
	\item When $p(x) = x - \lambda$ is linear with $\lambda \in K^{\times}$, $m_L^{\sigma}(x - \lambda)$ has a concrete interpretation in terms of norm equations and cyclic algebras (Section \ref{sect:linearcase}).
	\item For general $p(x) \in \sP_L$ with coefficients in $K$, we relate $m_L^{\sigma}(p(x))$ to the linear case above, but for a potentially different field corresponding to a root of $p(x)$ in $\overline{L}$ (Theorem \ref{thm:schurindexpx}).
	\item Finally, we show that the irreducible case when $t = 1$ suffices in general, by a lifting argument which shows that $m_L^{\sigma}(p(x)^t) = m_L^{\sigma}(p(x))$ for any $t \geq 1$ (Theorem \ref{thm:mpequalsmpt}).
\end{enumerate}

The main result, and most general statement, is Theorem \ref{thm:mainthm}.

\subsection{Reduction to Polynomials Over $K$}

\begin{defn}
	For $p(x) \in \sP_L$, let:
\begin{itemize}
		\item $\Gamma_{p(x)}$ denote the stabiliser of $p(x)$ in $\Gal(L/K)$.
		\item $K_{p(x)} \coloneqq L^{\Gamma_{p(x)}}$ denote the subfield of $L$ generated by $K$ and the coefficients of $p(x)$.
		\item $r(p) \coloneqq [\Gal(L/K) : \Gamma_{p(x)}] = [K_{p(x)} : K]$, so $r(p) \bZ$ is the preimage of $\Gamma_{p(x)}$ in $\bZ$.
	\end{itemize}
\end{defn}

The following is a special case of the general fact that to compute semilinear Schur indices, one can always ``factor out'' the Galois action. Here it says that to compute $m_L^{\sigma}(p(x)^t)$ in general, it suffices to understand the case when $p(x)$ has all its coefficients in $K$.

\begin{lemma}\label{lem:reductiontopolyoverK}
	For any $p(x) \in \sP_L$ and $t \geq 1$,
	\[
		m_L^{\sigma}(p(x)^t) = m_L^{\sigma^{r(p)}}(p(x)^t).
	\]
\end{lemma}

Here $m_L^{\sigma^{r(p)}}(p(x)^t)$ is defined as in Definition \ref{def:schurindexnotation} as the semilinear Schur index of $V_{p(x),t, s} \in \Ind_{L}(s \bZ)$ with respect to the action of $\bZ$ on $L$ where $1 \in \bZ$ acts as $\sigma^{r(p)}$, where $s \coloneqq r / r(p)$ is the order of $\sigma^{r(p)}$ in $\Aut(L)$. 

Equivalently, identifying $\bZ$ with $r(p)\bZ$, and thus the subgroup $s \bZ$ with $r \bZ$, $m_L^{\sigma^{r(p)}}(p(x)^t)$ is the semilinear Schur index of $V_{p(x),t, r} \in \Ind_{L}(r \bZ)$ with respect to the restricted action of $r(p)\bZ$ on $L$, where the generator $r(p) \in \bZ$ acts as $\sigma^{r(p)}$.

Either way, the key point is that the induced Galois action of $\Gal(L / K_{p(x)})$ on $\sP_L$ fixes $p(x)^t$.

\begin{proof}
	By definition, $m_L^{\sigma}(p(x)^t) = m^L_K(V_{p(x),t,r})$, where $V_{p(x),t,r}$ corresponds to $(p(x),t)$ under the bijection of Section \ref{sect:backgroundsemilinearreps}. This is compatible with the action of $\Gal(L/K)$, and therefore the stabiliser in $\Gal(L/K)$ of $V_{p(x),t,r}$ is the same as the stabiliser of $(p(x),t)$. This is equal to $\Gamma_{p(x)}$, the stabiliser of $p(x)$, and therefore by \cite[Prop.\ 7.3(7)]{TAYSchur} (in the notion of \emph{loc.~cit.})
	\[
	m_K^L(V_{p(x),t}) = m_{K_{p(x)}}^L(V_{p(x),t,r}),
	\]
	the second term of which is nothing but $m_L^{\sigma^{r(p)}}(p(x)^t)$ by the above discussion.
\end{proof}

\subsection{Matrix Description}

We can rephrase the condition that $m_L^{\sigma}(p(x)^t) \mid k$ purely in terms of matrices.

\begin{lemma}\label{lem:mdivkmatrixform}
    Suppose that $p(x) \in \sP_L$ has coefficients in $K$, and $t \geq 1$. Then for any $k \geq 1$, $m_L^{\sigma}(p(x)^t) \mid k$ if and only if there is some $A \in \GL_{knt}(L)$ with 
    \[
    N_{knt}^{\sigma}(A) = C(p^t)^{\oplus k}.
    \]
\end{lemma}

\begin{proof}
    To ease notation, set $m \coloneqq m_L^{\sigma}(p(x)^t)$. Suppose first there exists such an $A \in \GL_{knt}(L)$. Then the semilinear representation $W$ it defines, via $1 \in \bZ$ acting via $A$, satisfies $W|_{r\bZ} \cong (V_{p(x),t,r})^k$, and therefore $m \mid k$ by \cite[Prop.\ 9.3(4)]{TAYSchur}.
    
    Conversely, if $m \mid k$, then $k = a m$ for some $a \geq 1$. Let $A_0$ be the matrix of $1 \in \bZ$ acting on $W_{p(x),t, \sigma}$ with respect to some $L$-basis of $W_{p(x),t,\sigma}$. Now $N^{\sigma}_{mnt}(A_0)$ is the matrix of $r \in \bZ$ acting on $W_{p(x),t,\sigma}$, and therefore, as $(W_{p(x),t, \sigma})|_{r \bZ} \cong V_{p(x),t,r}^{\oplus m}$, there is some $P \in \GL_{mnt}(L)$ with 
    \[
    P N^{\sigma}_{mnt}(A_0) P^{-1} = C(p^t)^{\oplus m}.
    \]
    From the definition of $N^{\sigma}_{mnt}(-)$, we have that
    \[
        P N^{\sigma}_{mnt}(A_0) P^{-1} = N^{\sigma}_{mnt}(PA_0\sigma(P)^{-1}),
    \]
    and thus we may take $A$ to be the block sum of $a$ copies of $PA_0\sigma(P)^{-1}$.
\end{proof}

\begin{eg}
    For $p(x)$ as in Lemma \ref{lem:mdivkmatrixform}, it is always true that $m_L^{\sigma}(p(x)^t) \mid r$, by \cite[Prop.\ 9.3(6)]{TAYSchur}, and thus by Lemma \ref{lem:mdivkmatrixform} there must be some $A \in \GL_{rnt}(L)$ with $N_{rnt}^{\sigma}(A) = C(p^t)^{\oplus r}$. As a consistency check, we can find such an $A$ explicitly, taking
    \[
    A = \begin{pmatrix}
    0 & 0 &\cdots & 0 & C(p^t) \\
    I_{nt} & 0 & \cdots & 0 & 0 \\
    \vdots & \vdots & \ddots & \ddots & \vdots \\
    0 & 0 & \cdots & 0 & 0 \\
    0 & 0 & \cdots & I_{nt} & 0
\end{pmatrix},
\]
which satisfies $N_{rnt}^{\sigma}(A) = A^r = C(p^t)^{\oplus r}$.
\end{eg}

Finally, we record here a lemma that we shall use repeatedly in what follows.

\begin{lemma}\label{lem:commlemma}
	Suppose that $m \geq 1$, and $A \in M_m(L)$. If $N_m^{\sigma}(A) \in M_m(K)$, then $A N_m^{\sigma}(A) = N_m^{\sigma}(A) A$. The converse holds if $A$ is invertible.
\end{lemma}

\begin{proof}
	We may compute that
	\[
	A \sigma(N_m^{\sigma}(A)) = A \sigma(A) \cdots \sigma^{r-1}(A) A = N_m^{\sigma}(A)  A,
	\]
	from which the claim follows directly.
\end{proof}

\subsection{Linear Polynomials}\label{sect:linearcase}

For $k \geq 1$, we can rephrase the condition that $m_L^{\sigma}(x-\lambda) \mid k$ purely in terms of matrices using Lemma \ref{lem:mdivkmatrixform} (taking $n = 1$ and $t = 1$).

\begin{lemma}\label{lem:linearmatrixform}
    Suppose that $p(x) = x - \lambda$, where $\lambda \in K^{\times}$. Then for any $k \geq 1$, $m_L^{\sigma}(x-\lambda) \mid k$ if and only if there is some $A \in \GL_k(L)$ with 
    \[
    N_{k}^{\sigma}(A) = \lambda I_k.
    \]
    
    In particular, $m_L^{\sigma}(x - \lambda) = 1$ if and only if $\lambda \in N_{L/K}(L^{\times})$.
\end{lemma}

\begin{eg}
	For $\lambda \in K^{\times}$, $m_L^{\sigma}(x-\lambda) = 1$ whenever $L$ is finite by \cite[Cor.\ 9.4]{TAYSchur}, and so Lemma \ref{lem:linearmatrixform} recovers the classical fact that $N_{L/K} \colon L^{\times} \rightarrow K^{\times}$ is surjective when $L$ is a finite field.
\end{eg}

Lemma \ref{lem:linearmatrixform} provides the following necessary condition for $m_L^{\sigma}(x-\lambda)$ to divide $k$. We will see below (Remark \ref{rem:notsufficient}) that this is not sufficient in general, but turns out to be for certain fields.

\begin{cor}\label{cor:neccondition}
    Let $\lambda \in K^{\times}$. Then for any $k \geq 1$, if $m_L^{\sigma}(x-\lambda) \mid k$, then $\lambda^k \in N_{L/K}(L^{\times})$.
    
    Equivalently, $\lambda^{m_L^{\sigma}(x-\lambda)} = 1$ in $K^{\times} / N_{L/K}(L^{\times})$.
\end{cor}

\begin{proof}
    This follows directly from Lemma \ref{lem:linearmatrixform}, taking the determinant of $N_{k}^{\sigma}(A) = \lambda I_k$.
\end{proof}

We can also give a more theoretical description of $m_L^{\sigma}(x-\lambda)$ for $\lambda \in K^{\times}$. Recall that, as $L / K$ is Galois and $\lambda \in K$, we may form the cyclic algebra
\[
(L/K, \sigma, \lambda) \coloneqq L \oplus L \cdot u \oplus \cdots \oplus L \cdot u^{r-1}
\]
for some formal symbol $u$, which satisfies the relations
\[
u a = \sigma(a) u, \qquad u^r = \lambda,
\]
for $a \in L$. This is a central simple algebra over $K$ \cite[Thm.\ 14.6]{LAM}. 

\begin{defn}
    For $\lambda \in K^{\times}$, in order to ease notation, we write
    \begin{itemize}
    	\item $W_{\lambda} \coloneqq W_{(x-\lambda), 1, \sigma} \in \Irr_L^{\hspace{0.05em}  \sigma}(\bZ)$,
    	\item $V_{\lambda} \coloneqq V_{(x-\lambda), 1, r} \in \Irr_L(r\bZ)$.
    \end{itemize}
\end{defn}

\begin{lemma}\label{lem:linearbrauerequiv}
    For any $\lambda \in K^{\times}$, there is a $K$-algebra isomorphism
    \[
    (L/K, \sigma, \lambda) \cong M_{\frac{r}{m_L^{\sigma}(x-\lambda)}}(\End_{L \rtimes \bZ}(W_{\lambda})^{\opp}).
    \]
    In particular, $\End_{L \rtimes \bZ}(W_{\lambda})$ has centre $K$, and
    \[
    [\End_{L \rtimes \bZ}(W_{\lambda})] = - [(L/K, \sigma, \lambda)]
    \]
    in $\Br(K)$.
\end{lemma}

\begin{proof}
    As $(L/K, \sigma, \lambda)$ is a central simple algebra over $K$, $(L/K, \sigma, \lambda) \cong M_n(D^{\opp})$ as $K$-algebras, for some unique $K$-division algebra $D$. The $K$-linear category of $M_n(D^{\opp})$-modules is semisimple, with a unique simple object $U$, and the $K$-division algebra $D$ is recovered as the endomorphism ring of $U$. It therefore suffices to give $W_{\lambda}$ the structure of a $(L/K,  \sigma, \lambda)$-module such that
    \[
    \End_{(L/K,  \sigma, \lambda)}(W_{\lambda}) = \End_{L \rtimes \bZ}(W_{\lambda}).
    \]
    
    There is a canonical identification
    \[
    L_\sigma[x, x^{-1}] \xrightarrow{\sim} L \rtimes_{\sigma} \bZ, \qquad x \mapsto 1 \in \bZ,
    \]
    which restricts to an identification
    \[
    L[x^r, x^{-r}] \xrightarrow{\sim} L[r \bZ].
    \]
    The resulting $L_{\sigma}[x, x^{-1}]$-module structure on $W_{\lambda}$ is killed by $x^r - \lambda$, because $W_{\lambda}$ restricts to $V_{\lambda}^{\oplus m_L^{\sigma}(x - \lambda)}$ as an $L[r \bZ]$-module, and $V_{\lambda}$ is killed by $x^r - \lambda$.

    This induces a $(L/K, \sigma,  \lambda)$-module structure on $W_{\lambda}$, via the isomorphism
    \[
    (L/K,  \sigma, \lambda) \xrightarrow{\sim} L_{\sigma}[x,x^{-1}] / (x^r - \lambda),
    \]
    which has the property that $\End_{(L/K,  \sigma,  \lambda)}(W_{\lambda}) = \End_{L \rtimes_{\sigma} \bZ}(W_{\lambda})$, as required.
\end{proof}

Recall that for a class $[A] \in \Br(K)$, the \emph{exponent} (or \emph{period}) $\exp [A]$ is the order of $[A]$ in $\Br(K)$, and the \emph{index} $\ind [A]$ is $\sqrt{\dim_K D}$ for the unique central division algebra $D$ over $K$ with $[A] = [D]$.

\begin{cor}\label{cor:cyclicindex}
    Suppose that $\lambda \in K$. Then 
    \[
    m_L^{\sigma}(x - \lambda) = \ind [(L/K, \sigma, \lambda)].
    \]
\end{cor}

\begin{proof}
From the isomorphism of Lemma \ref{lem:linearbrauerequiv}
    \[
    (L/K,  \sigma, \lambda) \cong M_{\frac{r}{m_L^{\sigma}(x-\lambda)}}(\End_{L \rtimes \bZ}(W_{\lambda})^{\opp})
    \]
    we have that $\End_{L \rtimes \bZ}(W_{\lambda})$ has centre $K$, and therefore $m_L^{\sigma}(x-\lambda)$ is the index of $\End_{L \rtimes \bZ}(W_{\lambda})$ by \cite[Prop.\ 8.29(2)]{TAYSchur}, and thus from the isomorphism also the index of $(L/K,  \sigma, \lambda)$.
\end{proof}

We can use Corollary \ref{cor:cyclicindex} to recover and reinterpret many classical properties of cyclic algebras from properties of semilinear Schur indices. The following are Theorems 14.7, 14.8 and 14.9 respectively of \cite{LAM}.

\begin{cor}\label{cor:cyclicalgproperties}
    Suppose that $\lambda \in K$. Then:
    \begin{enumerate}
        \item $(L/K,  \sigma, \lambda) \cong M_r(K)$ if and only if $\lambda \in N_{L/K}(L^{\times})$.
        \item If $r$ is prime, then $(L/K, \sigma,  \lambda)$ is a division algebra if and only if $\lambda \not\in N_{L/K}(L^{\times})$.
        \item If $\lambda, \lambda^2, ... ,\lambda^{r-1} \not\in N_{L/K}(L^{\times})$, then $(L/K,  \sigma, \lambda)$ is a division algebra.
    \end{enumerate}
\end{cor}

\begin{proof}
        For the first point, note that $(L/K,  \sigma, \lambda) \cong M_r(K)$ if and only if $(L/K,  \sigma, \lambda)$ has index $1$, and therefore by Corollary \ref{cor:cyclicindex} if and only if $m_L^{\sigma}(x - \lambda) = 1$, or equivalently that $\lambda \in N_{L/K}(L^{\times})$ by Lemma \ref{lem:linearmatrixform}. For the second and third points, note that $(L/K,  \sigma, \lambda)$ is a division algebra if and only if $(L/K,  \sigma, \lambda)$ has index $r$, or equivalently that $m_L^{\sigma}(x - \lambda) = r$ by Corollary \ref{cor:cyclicindex}. Furthermore, as is generally true, we have that $m_L^{\sigma}(x - \lambda) \mid r$ by \cite[Prop.\ 9.3(6)]{TAYSchur}. If $r$ is prime, the result thus follows from the first point. If instead $\lambda, \lambda^2, ... ,\lambda^{r-1} \not\in N_{L/K}(L^{\times})$, then $m_L^{\sigma}(x - \lambda) \nmid k$ for any $1 \leq k < r$ by Corollary \ref{cor:neccondition}, and thus to only option is that $m_L^{\sigma}(x - \lambda) = r$.
\end{proof}

The third point of Corollary \ref{cor:cyclicalgproperties} is sometimes known as Wedderburn's Theorem, and provides a sufficient condition for $(L/K,  \sigma, \lambda)$ to be a division algebra. This condition is not necessary in general. We can use Lemma \ref{lem:linearbrauerequiv} to state an equivalent condition for $(L/K,  \sigma, \lambda)$ to be a division algebra, which recovers Wedderburn's Theorem upon taking the determinant.

\begin{cor}
    Suppose that $\lambda \in K^{\times}$. Then for any $k \geq 1$, $\ind(L/K,  \sigma, \lambda) \mid k$ if and only if $\lambda I_k \in N_k^{\sigma}(\GL_k(L))$. In particular, $(L/K,  \sigma, \lambda)$ is a division algebra if and only if
    \[
    \lambda I_k \not\in N_k^{\sigma}(\GL_k(L))
    \]
    for $k = 1, ... , r-1$.
\end{cor}

\begin{proof}
    The first claim follows from Lemma \ref{lem:linearmatrixform} and Lemma \ref{lem:linearbrauerequiv}. The second claim follows from the first, as $(L/K, \sigma,  \lambda)$ is a division algebra if and only if $\ind(L/K,  \sigma, \lambda) = r$, and $m_L^{\sigma}(x - \lambda) = \ind(L/K,  \sigma, \lambda)$ by Corollary \ref{cor:cyclicindex}.
\end{proof}

Lemma \ref{lem:linearbrauerequiv} also allows us to use results about cyclic algebras to make semilinear Schur indices more computable over certain nice fields.

\begin{cor}\label{cor:schurindexisorder}
    Suppose that $K$ is a field such that the period and the index of elements of $\Br(K)$ coincide. Then $m_L^{\sigma}(x- \lambda)$ is the order of $\lambda$ in $K^{\times} / N_{L / K}(L^{\times})$. Equivalently, for any $k \geq 1$, $m_L^{\sigma}(x - \lambda) \mid k$ if and only if $\lambda^k \in N_{L / K}(L^{\times})$.
\end{cor}

\begin{proof}
    Because $L / K$ is cyclic, the map
    \[
        K^{\times} \rightarrow \Br(L/K), \qquad \lambda \mapsto [(L/K,  \sigma, \lambda)]
    \]
    is surjective with kernel $N_{L/K}(L^{\times})$ \cite[\S 15.1, Prop.\ b]{PIE}. In particular, we see that the order of $\lambda$ in $K^{\times} / N_{L / K}(L^{\times})$ agrees with $\exp([A])$, and so by assumption with $\ind([(L/K, \sigma,  \lambda)])$, which is simply $m_L^{\sigma}(x- \lambda)$ by Corollary \ref{cor:cyclicindex}.
\end{proof}

\begin{remark}\label{rem:localandglobBrauerGroups}
Fields which satisfy the assumptions of Corollary \ref{cor:schurindexisorder} include local fields \cite[Rem.\ IV.4.4(b)]{milneCFT} and global fields \cite[Thm.\ VIII.2.6]{milneCFT}, and so for such fields the computation of the semilinear Schur indices $m_L^{\sigma}(x-\lambda)$ becomes fully algorithmic by computing in $K^{\times} / N_{L/K}(L^{\times})$.
\end{remark}

\begin{remark}\label{rem:notsufficient}
Corollary \ref{cor:schurindexisorder} allows us to see that the converse of Corollary \ref{cor:neccondition} is not true in general. Indeed, any $\lambda$ for which $\exp[(L/K,  \sigma, \lambda)] \neq \ind[(L/K,  \sigma, \lambda)]$ will have $\lambda$ of order strictly less than $m_L^{\sigma}(x-\lambda)$ in $K^{\times} / N_{L/K}(L^{\times})$.
\end{remark}

\subsection{Higher Degree Polynomials}\label{sect:higherdegree}

We suppose in this section that $p(x) \in \sP_L$, with coefficients in $K$. Let $n \geq 1$ denote the degree of $p(x)$, and let $\alpha$ be a root of $p(x)$ in $\overline{L}$.

\begin{defn}\label{def:taualpha}
We set $\tau_{\alpha}$ to be the automorphism of $L(\alpha)$, which is $\sigma$ on $L$ and trivial on $\alpha$.
\end{defn}

This is well-defined because $p(x)$ has coefficients in $K$, so the automorphism $\sigma$ of $L[x]$ descends to an automorphism of $L[x] / p(x) = L(\alpha)$. The order of $\tau_{\alpha}$ is $r$, as the restriction $\sigma$ has order $r$ on $L$, and thus $L(\alpha) / K(\alpha)$ is Galois, with Galois group $\langle \tau_{\alpha} \rangle$ of order $r$. There is thus an action of $\bZ$ on $L(\alpha)$, where $1 \in \bZ$ acts as $\tau_{\alpha}$, and we may consider semilinear representations of $\bZ$ over $L(\alpha)$.

The aim of this section is to prove the following, which shows that the computation of semilinear Schur indices is reduced to the case of linear polynomials, over a potentially different field. 

\begin{thm}\label{thm:schurindexpx}
	Suppose that $p(x) \in \sP_L$, and has coefficients in $K$. Then
	\[
	m_L^{\sigma}(p(x)) = m_{L(\alpha)}^{\tau_{\alpha}}(x - \alpha).
	\]
\end{thm}

Towards the proof of Theorem \ref{thm:schurindexpx}, note that there is an $L$-algebra homomorphism
\[
M(-): L(\alpha) \rightarrow M_n(L),
\]
defined by sending any $x \in L(\alpha)$ to the matrix of left multiplication by $x$ on $L(\alpha)$ with respect to the ordered $L$-basis $\{1, \alpha, \cdots , \alpha^{n-1}\}$. Note that $M(\alpha) = C(p)$, the companion matrix of $p(x)$.

\begin{lemma}\label{lem:imageoffield1dim}
   The $L$-algebra homomorphism 
   \[
   M(-): L(\alpha) \rightarrow M_n(L)
   \]
   is injective, the image consists of those matrices which commute with $M(\alpha)$, and for any $x \in L(\alpha)$,
   \[
   M(\tau_{\alpha}(x)) = \sigma(M(x)).
   \]
\end{lemma}

\begin{proof}
   We can consider $M(-)$ as the composition
   \[
   L(\alpha) \xrightarrow{m_{(-)}} \End_L(L(\alpha)) \xrightarrow{\sim} M_n(L),
   \]
   where for $x \in L$, $m_x \in \End_L(L(\alpha))$ is left-multiplication by $x$, and so $M(-)$ is injective. The description of the image follow from Lemma \ref{lem:centralisercyclicmodule} because $M(\alpha) = C(p)$, so the image is exactly $L[C(p)]$. For the final point, note that $\sigma(M(x))$ is in the image of $M(-)$, as $\sigma(M(\alpha)) = M(\alpha)$ so
\[
\sigma(M(x)) M(\alpha) = \sigma(M(x)M(\alpha)) = \sigma(M(\alpha)M(x)) = M(\alpha) \sigma(M(x)).
\]
There is thus some $y \in L(\alpha)$ with $M(y) = \sigma(M(x))$. We can recover $y$ as
\[
y = y_0 + y_1 \alpha + \cdots + y_{n-1}\alpha^{n-1},
\]
where $(y_0, \cdots, y_{n-1})^T$ is the first column of $\sigma(M(x))$. Now $M(x)$ has first column $(x_0, \cdots, x_{n-1})^T$, where
\[
x = x_0 + x_1\alpha + \cdots + x_{n-1}\alpha^{n-1},
\]
so $\sigma(M(x))$ has first column $(\sigma(x_0), \cdots, \sigma(x_{n-1}))^T$. We may thus compute
\[
y = \sigma(x_0) + \sigma(x_1) \alpha + \cdots + \sigma(y_{n-1})\alpha^{n-1} = \sigma(x_0 + x_1\alpha + \cdots + x_{n-1}\alpha^{n-1}) = \tau_{\alpha}(x),
\]
and therefore $M(\tau_{\alpha}(x)) = \sigma(M(x))$.
\end{proof}

We may extend $M(-)$ to matrices over $L(\alpha)$ as follows.

\begin{defn}
   For $m \geq 1$, we define
   \[
   M(-) \colon M_m(L(\alpha)) \rightarrow M_{mn}(L),
   \]
   by setting, for $A = (a_{ij})_{i,j} \in M_m(L(\alpha))$,
   \[
       M(A) \coloneqq \begin{pmatrix}
   M(a_{11}) & \cdots & M(a_{1m}) \\
   \vdots & \ddots & \vdots \\
   M(a_{m1})  & \cdots & M(a_{mm})
\end{pmatrix}.
   \]
\end{defn}

It is direct to see that $M(-) \colon M_m(L(\alpha)) \rightarrow M_{mn}(L)$ is an injective $L$-algebra homomorphism, which satisfies $M(\tau_{\alpha}(A)) = \sigma(M(A))$ for any $A \in M_m(L(\alpha))$.

\begin{lemma}\label{lem:bijsolsetst1}
   For any $k \geq 1$, $M(-)$ induces a bijection
   \[
   M(-) \colon \left\{A \in M_k(L(\alpha)) \mid N_k^{\tau_{\alpha}}(A) = \alpha I_k \right\} \xrightarrow{\sim} \left\{ B \in M_{kn}(L) \mid N_{kn}^{\sigma}(B) = C(p)^{\oplus k} \right\}.
   \]
\end{lemma}

\begin{proof}
   By the above properties of $M(-)$, then for any $A \in M_k(L(\alpha))$ with $N_k^{\tau_{\alpha}}(A) = \alpha I_k$,
   \begin{align*}
   N_{kn}^{\sigma}(M(A)) &= M(A)\sigma(M(A)) \cdots \sigma^{r-1}(M(A)), \\
   &= M(A)M(\tau_{\alpha}(A)) \cdots M(\tau_{\alpha}^{r-1}(A)), \\
   &= M(N_{k}^{\tau_{\alpha}}(A)), \\
   &= M(\alpha I_k), \\
   &= C(p)^{\oplus k}.
   \end{align*}
   This shows that the map is well-defined and injective. 

   For surjectivity, suppose that $B \in M_{kn}(L)$ with $N_{kn}^{\sigma}(B) = C(p)^{\oplus k}$. By Lemma \ref{lem:commlemma}, $B$ commutes with $C(p)^{\oplus k}$, as $C(p)^{\oplus k}$ is valued over $K$. If we consider $B$ as a block matrix, so an element of $M_k(M_n(L))$, we see that each entry commutes with $C(p) = M(\alpha)$, and therefore by Lemma \ref{lem:imageoffield1dim} $B = M(A)$ for some $A \in M_k(L(\alpha))$. 

   Similarly, to above, the condition on $B$ shows that
   \[
   M(N_k^{\tau_{\alpha}}(A)) = N_{kn}(M(A)) = N_{kn}(B) = C(p)^{\oplus k} = M(\alpha I_k),
   \]
   and therefore $N_k^{\tau_{\alpha}}(A) = \alpha I_k$ from the injectivity of $M(-)$.
\end{proof}

We can use this to deduce Theorem \ref{thm:schurindexpx}.

\begin{proof}[Proof of Theorem \ref{thm:schurindexpx}]
   We have that the first (resp.\ the second) set appearing in the statement of Lemma \ref{lem:bijsolsetst1} is non-empty if and only if $m_{L(\alpha)}^{\tau_{\alpha}}(x - \alpha) \mid k$ (resp.\  $m_L^{\sigma}(p(x)) \mid k$) by Lemma \ref{lem:mdivkmatrixform}.
   
   Therefore, for any $k \geq 1$, $m_{L(\alpha)}^{\tau_{\alpha}}(x - \alpha) \mid k$ if and only if $m_L^{\sigma}(p(x)) \mid k$ by Lemma \ref{lem:bijsolsetst1}, and thus $m_L^{\sigma}(p(x)) = m_{L(\alpha)}^{\tau_{\alpha}}(x - \alpha)$.
\end{proof}

\subsection{Indecomposable Representations}

The aim of this section is to extend Theorem \ref{thm:schurindexpx} from irreducible to indecomposable representations, by a lifting argument.

\begin{thm}\label{thm:mpequalsmpt}
Suppose that $p(x) \in \sP_L$ has coefficients in $K$. Then for any $t \geq 1$,
\[
m_L^{\sigma}(p(x)) = m_L^{\sigma}(p(x)^t).
\]
\end{thm}

\begin{proof}
	It suffices to show that for any $k \geq 1$, $m_L^{\sigma}(p(x)) \mid k$ if and only if $m_L^{\sigma}(p(x)^t) \mid k$, and so let us fix $k \geq 1$.
	
	Let, for $t \geq 1$, $S_t \coloneqq \fC(C(p^t))$, the centraliser of $C(p^t)$ in $M_{nt}(L)$, which is of interest as for any $t \geq 1$, $m_L^{\sigma}(p(x)^t) \mid k$ if and only if there is some $A \in M_{k}(S_t)$ with $N_{knt}^{\sigma}(A) =  C(p^t)^{\oplus k}$, by Lemma \ref{lem:mdivkmatrixform} and Lemma \ref{lem:commlemma}. Here we view $M_k(S_t) \subset M_{knt}(L)$ via the natural inclusion $S_t \subset M_{nt}(L)$ and identification $M_{k}(M_{nt}(L)) \xrightarrow{\sim} M_{knt}(L)$.

	We can describe $S_t$ more explicitly. For $t \geq 1$, set $R_t \coloneqq L[x] / \fp^t$, where $\fp = (p(x))$ is the principal ideal of $L[x]$ generated by $p(x)$. From Lemma \ref{lem:centralisercyclicmodule}, the map
	\[
	R_t \rightarrow S_t, \qquad x \mapsto C(p^t)
	\]
	is an $L$-algebra isomorphism, and thus induces an $L$-algebra isomorphism
	\[
		\phi_t \colon M_k(R_t) \xrightarrow{\sim} M_k(S_t).
	\]

	The natural action of $\Gal(L/K)$ on $L[x]$ induces an action on $R_t$ because $p(x)$ is fixed by $\Gal(L/K)$, and thus an action on $\Gal(L/K)$ on $M_k(R_t)$ entrywise. With respect to this action $\phi_t$ is $L \rtimes \Gal(L/K)$-linear, and so setting
	\[
	\sN_t \colon M_k(R_t) \rightarrow M_k(R_t), \qquad \sN_t(A) \coloneqq A \sigma(A) \cdots \sigma^{r-1}(A),
	\]
	we have that $\phi_t \circ \sN_t = N_{knt}^{\sigma} \circ \phi_t$. In particular, for any $t \geq 1$, we have rephrased the condition: $m_L^{\sigma}(p(x)^t) \mid k$ if and only if there is some $A_t \in M_{k}(R_t)$ with $\sN_{t}(A_t) = x I_k$.

	For any $t \geq 1$, there are natural quotient maps
	\[
	q_t \colon M_k(R_{t+1}) \rightarrow M_k(R_t)
	\]
	induced by the natural quotient map $R_{t+1} \twoheadrightarrow R_t$, which satisfy $\sN_t \circ q_t = q_t \circ \sN_{t+1}$. 

	In particular, if $m_L^{\sigma}(p(x)^t) \mid k$, then taking $A_t \in M_{k}(R_t)$ with $\sN_{t}(A_t) = x I_k$, $\sN_1(q_{t,1}(A_t))) = x I_k$, so $m_L^{\sigma}(p(x)) \mid k$, where
	\[
		q_{t,1} \coloneqq q_1 \circ \cdots \circ q_{t-1} \colon M_k(R_t) \twoheadrightarrow M_k(R_1).
	\]

	In the other direction, suppose that $m_L^{\sigma}(p(x)) \mid k$. We show that $m_L^{\sigma}(p(x)^t) \mid k$ by induction on $t \geq 1$. When $t = 1$ there is nothing to prove, so suppose that $m_L^{\sigma}(p(x)^t) \mid k$ for some $t \geq 1$, and let $A_t \in M_{k}(R_t)$ with $\sN_{t}(A_t) = x I_k$. Pick any lift $A \in M_k(R_{t+1})$ of $A_t$. We want to modify $A$, to obtain some $A_{t+1}$ with $\sN_{t+1}(A_{t+1}) = x I_k$. Towards this, let
	\[
		E \coloneqq \sN_{t+1}(A) - x I_k \in J
	\]
	denote the error. 
	
	Because $p(x)$ is coprime to $x$, $x I_k \in \GL_k(R_t)$, and thus $A_t \in \GL_k(R_t)$ too. The kernel $J = M_k(\fp^t / \fp^{t+1})$ of $q_t \colon M_k(R_{t+1}) \rightarrow M_k(R_t)$ is nilpotent, as $J^2 = 0$, and thus $A \in \GL_k(R_{t+1})$.
	
	Because $A$ is invertible, we may then define,
	\[
	\rho \colon J \rightarrow J, \qquad \rho(M) \coloneqq A \sigma(M) A^{-1},
	\]
	and 
	\[
	\Tr_{\rho} \colon J \rightarrow J, \qquad \Tr_{\rho} \coloneqq M + \rho(M) + \cdots + \rho^{r-1}(M).
	\]

	Using this, for a yet to be determined $D \in J$, we may compute the norm of our candidate perturbation $A_{t+1} \coloneqq (1+D)A$ of $A$:
	\begin{align*}
		\sN_{t+1}(A_{t+1}) &= (1+D)A \sigma(1+D) \sigma(A) \cdots \sigma^{r-1}(1+D) \sigma^{r-1}(A), \\
			&= (1+D)\rho(1+D) A\sigma(A)\sigma^2(1+D) \sigma^2(A) \cdots \sigma^{r-1}(1+D) \sigma^{r-1}(A), \\
			&= \sN_{t+1}(1+D) \sN_{t+1}(A), \\
			&= (1 + \Tr_{\rho}(D)) \sN_{t+1}(A), \\
			&= (1 + \Tr_{\rho}(D)) (x I_k + E), \\
			&= x I_k + E + x \Tr_{\rho}(D),
	\end{align*}
	where we have used that $\sN_{t+1}(1+D) = 1 + \Tr_{\rho}(D)$ and $\Tr_{\rho}(D)E = 0$ because $J^2 = 0$. We want that $\sN_{t+1}(A_{t+1}) = x I_k$, and therefore want to choose $D \in J$ with 
	\[
		\Tr_{\rho}(D) = - x^{-1} E,
	\]
	where $x^{-1}$ is the inverse of $x$ in $R_{t+1}$. The key property that allows us to do this is that $-x^{-1} E$ is fixed by $\rho$, which holds because $\rho(-x^{-1} E) = -x^{-1} \rho(E)$ and
	\[
		\rho(E) = A \sigma(E) A^{-1} = A(\sigma(\sN_{t+1}(A)) - x I_k)A^{-1} = A \sigma(A) \cdots \sigma^{r-1}(A) \sigma^r(A) A^{-1} - x I_k = E.
	\]
		Indeed, for any $M \in J$ with $\rho(M) = M$, taking $u \in L$ be any element with $\Tr_{L/K}(u) = 1$ (which is possible because $L / K$ is separable), we have that
	\begin{align*}
		\Tr_{\rho}(u M) &= u M + \rho(u M) + \cdots + \rho^{r-1}(uM) = (u + \sigma(u) + \cdots + \sigma^{r-1}(u))M = M.
	\end{align*}
	Explicitly, for $M := - x^{-1}E$, we may find the required $D \coloneqq -u x^{-1} E$, so $A_{t+1} = (1 - ux^{-1} E)A$.
\end{proof}

\subsection{General Statement}

For $p(x) \in \sP_L$, recall that $K_{p(x)}$ denotes the subfield of $L$ generated by $K$ and the coefficients of $p(x)$. We first generalise the definition of $\tau_{\alpha}$ of Section \ref{sect:higherdegree} to the case when $p(x)$ doesn't have coefficients in $K$. 

\begin{defn}
	We write $\tau_{\alpha}$ for the automorphism of $L(\alpha)$ which is $\sigma^{r(p)}$ on $L$, trivial on $\alpha$.
\end{defn}

This is simply Definition \ref{def:taualpha} for the extension $L/K_{p(x)}$, which is Galois, has Galois group generated by $\sigma^{r(p)}$, and $p(x)$ has coefficients in $K_{p(x)}$. In particular, $\tau_{\alpha}$ is well-defined, and we have that $L(\alpha) / K_{p(x)}(\alpha)$ is Galois, with $\langle \tau_{\alpha} \rangle = \Gal(L(\alpha) / K_{p(x)}(\alpha))$ of order $r / r(p)$.

We first note the following notational simplification.

\begin{lemma}\label{lem:dontneedpx}
	The subfields $K_{p(x)}(\alpha)$ and $K(\alpha)$ of $L(\alpha)$ are equal.
\end{lemma}

\begin{proof}
	For any $\gamma \in \Gal(L(\alpha) / K(\alpha))$, $\gamma(p)(\alpha) = \gamma(p(\alpha)) = 0$ because $\gamma$ fixes $\alpha$. Now $\gamma$ preserves the Galois subfield $L$, hence $\gamma(p)$ is a monic irreducible polynomial over $L$ with $\gamma(p)(\alpha) = 0$, thus $\gamma(p) = p$. This being true for any $\gamma$, all coefficients of $p(x)$ are thus contained in $K(\alpha)$.
\end{proof}

We can summarise the results of this section in the following theorem.

\begin{thm}\label{thm:mainthm}
	Suppose that $p(x) \in \sP_L$ and $t \geq 1$. Let $\alpha \in \overline{L}$ be a root of $p(x)$. Then 
	\[
	m_L^{\sigma}(p(x)^t) = m_{L(\alpha)}^{\tau_{\alpha}}(x - \alpha).
	\]
	If additionally the period equals the index in $\Br(K(\alpha))$, then $m_L^{\sigma}(p(x)^t)$ is the order of $\alpha$ in
	\[
		K(\alpha)^{\times} / N_{L(\alpha) / K(\alpha)}(L(\alpha)^{\times}).
	\]

	In particular, this holds if $L$ is a local or global field.
\end{thm}

\begin{proof}
From Lemma \ref{lem:reductiontopolyoverK}, Theorem \ref{thm:mpequalsmpt} and Theorem \ref{thm:schurindexpx} respectively, there are equalities
\[
m_L^{\sigma}(p(x)^t) = m_L^{\sigma^{r(p)}}(p(x)^t) = m_L^{\sigma^{r(p)}}(p(x)) = m_{L(\alpha)}^{\tau_{\alpha}}(x - \alpha).
\]
The second claim follows from Corollary \ref{cor:schurindexisorder}, and the final claim from Remark \ref{rem:localandglobBrauerGroups}.
\end{proof}

\subsection{Examples}

\subsubsection{The Extension $\bC / \bR$}\label{eg:CoverR}

Suppose that $L / K$ is $\bC / \bR$, so $\sigma = \overline{(-)}$. Then $\sP_{\bC} = \{x - \lambda \mid \lambda \in \bC^{\times}\}$ is canonically identified with $\bC^{\times}$. In Theorem \ref{thm:mainthm}, $\bC(\alpha) = \bC$ and $\alpha = \lambda$, hence $m_{\bC}^{\sigma}(x - \lambda)$ is simply the order of $\lambda$ in $\bR(\lambda)^{\times} / N_{\bC / \bR(\lambda)}(\bC^{\times})$. There are two cases to consider:
\begin{enumerate}
		\item $\lambda \not\in \bR$, so $\bR(\lambda) = \bC$, so
		\[
		\bR(\lambda)^{\times} / N_{\bC / \bR(\lambda)}(\bC^{\times}) = \bC^{\times} / N_{\bC / \bC}(\bC^{\times}) = 1
		\]
		and hence $m_{\bC}^{\sigma}(x - \lambda) = 1$.
		\item If $\lambda \in \bR$, then $\bR(\lambda) = \bR$, so 
		\[
		\bR(\lambda)^{\times} / N_{\bC / \bR(\lambda)}(\bC^{\times}) = \bR^{\times} / N_{\bC / \bR}(\bC^{\times}) = \bR^{\times} / \bR_{> 0}
		\]
		and therefore
        \[
            m_{\bC}^{\sigma}(x - \lambda) = \twopartdef{1}{\lambda > 0}{2}{\lambda < 0}.
        \]
	\end{enumerate}
Thus in general, for $\lambda \in \bC^{\times}$ and $t \geq 1$,
        \[
        m_{\bC}^{\sigma}(x - \lambda) = \twopartdef{1}{\lambda \not\in \bR_{< 0}}{2}{\lambda \in \bR_{< 0}}.
        \]

\subsubsection{Finite Fields}

If $L$ is finite, then $m_L^{\sigma}(p(x)^t) = 1$ for any $p(x) \in \sP_L$ and $t \geq 1$. This follows directly from Theorem \ref{thm:mainthm}, as the norm map is surjective for an extension of finite fields.

\subsection{Endomorphism Rings}\label{sect:endomorphimrings}

We can also use the methods described above to describe the endomorphims rings of the irreducible representation $W_{p(x), 1}$ for $p(x) \in \sP_L$. We will use this in Section \ref{sect:localfields} to describe the endomorphism rings of irreducible isocrystals over finite fields.

To ease notation, for $W \in \Rep_L^{\sigma}(\bZ)$ we write $\End_{L}^{\sigma}(W) \coloneqq \End_{L \rtimes_{\sigma} \bZ}(W)$ for the endomorphism ring of $W$ in $\Rep_L^{\sigma}(\bZ)$, which is naturally a $K$-algebra. We also set $W_{p(x)} \coloneqq W_{p(x), 1}$ for $p(x) \in \sP_L$.

We first describe how $\End_{L}^{\sigma}(W_{p(x)})$ canonically has the structure of a $K(\alpha)$-algebra, which is not immediately clear.

\begin{lemma}\label{lem:generalcentre}
	Suppose that $W \in \Rep_L^{\sigma}(\bZ)$. Then $Z(\End_L^{\sigma}(W)) = K[\rho_W(r)]$, where $K[\rho_W(r)]$ is the $K$-subalgebra of $\End_L^{\sigma}(W)$ generated by $\rho_W(r)$.
\end{lemma}

\begin{proof}
	We have that $K[\rho_W(r)] \subset Z(\End_L^{\sigma}(W))$, and so it suffices to show that the induced map $L \otimes_K K[\rho_W(r)] \rightarrow L \otimes_K Z(\End_L^{\sigma}(W))$ is surjective. From \cite[Thm.\ A]{TAYSchur} and \cite[Prop. IV.2.3]{milneCFT} there is an identification
	\begin{equation}\label{eqn:basechangeiso}
	L \otimes_K Z(\End_L^{\sigma}(W)) \xrightarrow{\sim} Z(\End_{L[r\bZ]}(W|_{r \bZ})),
	\end{equation}
	and the image of $L \otimes_K K[\rho_W(r)]$ is the subalgebra $L[\rho_W(r)] \subset Z(\End_{L[r\bZ]}(W|_{r \bZ}))$. That this is an equality follows from the general fact that for $M$ a finite-dimensional $L[x]$-module, the natural map $L[x] \rightarrow Z(\End_{L[x]}(M))$ is surjective. Indeed, suppose that
	\[
	M = L[x] / f_1 \oplus \cdots \oplus L[x] / f_m
	\]
	with $f_1 \mid f_2 \mid \cdots \mid f_m$. For $\phi \in Z(\End_{L[x]}(M))$, $\phi$ must commute with the projection maps to each factor, hence $\phi = (m_{q_1}, ... ,m_{q_m})$ is given by left-multiplication maps on each factor for some $q_1, ... ,q_m \in L[x]$. Then commutativity with each reduction map $L[x]/f_m \rightarrow L[x] / f_i$ shows that $q_i \equiv q_m \mod (f_i)$ for all $1 \leq i \leq m$, hence $\phi$ is simply left-multiplication by $q_m \in L[x]$.
\end{proof}

When $W$ is irreducible, $\End_L^{\sigma}(W)$ is a division ring and so $K[\rho_W(r)]$ is a field. We can explicitly describe this as follows. Throughout, we fix $\alpha \in \overline{L}$, a root of $p(x)$.

\begin{lemma}\label{lem:centredesc}
	Suppose that $p(x) \in \sP_L$, and $W \in \Rep_{L}^{\sigma}(\bZ)$ with $W \cong W_{p(x)}$. Then there is a canonical $K$-algebra isomorphism
	\[
	K(\alpha) \xrightarrow{\sim} Z(\End_L^{\sigma}(W)), \qquad \alpha \mapsto \rho_W(r).
	\]
	If $W_1 \cong W_{p(x)} \cong W_2$, then for any isomorphism $f \colon W_1 \xrightarrow{\sim} W_2$, the induced map
	\[
	f_* \colon \End_L^{\sigma}(W_1) \rightarrow \End_L^{\sigma}(W_2), \qquad \phi \mapsto f \circ \phi \circ f^{-1}
	\]
	is an isomorphism of $K(\alpha)$-algebras with respect to these $K(\alpha)$-algebra structures.
\end{lemma}

\begin{proof}
	From Lemma \ref{lem:generalcentre}, $Z(\End_L^{\sigma}(W)) = K[\rho_W(r)]$ is a field. Let $g(y) \in K[y]$ be the minimal polynomial of $\rho_W(r)$ over $K$. As $W \cong W_{p(x)}$, the from Theorem \ref{thm:descofrestriction} $W|_{r \bZ}$ has $V_{p(x), 1, r} = L[x] / p(x)$ as a direct summand. The action of $\rho_W(r)$ on $L[x] / p(x)$ is by multiplication by $x$, hence as $g(\rho_W(r)) = 0$, $g(x)$ is zero in $L[x] / p(x)$, and so $p(x) \mid g(x)$. This shows that $g(\alpha) = 0$, hence $g(y)$ is the minimal polynomial of $\alpha$ over $K$, and thus the required map is an isomorphism.

	The second claim follows directly, using that $f \circ \rho_{W_1}(r) = \rho_{W_2}(r) \circ f$.
\end{proof}

We now describe $\End_{L}^{\sigma}(W_{p(x)})$ with respect to this $K(\alpha)$-algebra structure. We first need to relate $\End_{L}^{\sigma}(W_{p(x)})$ to an endomorphism ring of an associated representation of the subgroup $r(p) \bZ$.

Let $p(x) \in \sP_L$ and $W \in \Rep_L^{\sigma}(\bZ)$ with $W \cong W_{p(x)}$. Recall that $r(p) = [K_{p(x)} : K] \mid r$, and so $r \bZ \leq r(p) \bZ \leq \bZ$. Recall that $\Gamma_{p(x)}$ is the stabiliser of $p(x)$ in $\Gamma = \Gal(L/K)$, so
\[
W|_{r \bZ} \cong \left( \bigoplus_{\gamma \in \Gamma / \Gamma_{p(x)}} \gamma * V_{p(x), 1, r} \right)^{\oplus m_L^{\sigma}(p(x))}.
\]
From \cite[Cor.\ 8.27]{TAYSchur}, 
\[
W|_{r(p) \bZ} \cong \bigoplus_{\gamma \in \Gamma / \Gamma_{p(x)}} U_{\gamma},
\]
where each $U_{\gamma}$ is an irreducible semilinear representation of $r(p)\bZ$ over $L$ (where $r(p)$ acts by $\sigma^{r(p)}$ on $L$), which are pairwise non-isomorphic and $U_{\gamma}|_{r \bZ} \cong (\gamma * V_{p(x), 1, r})^{\oplus m_L^{\sigma}(p(x))}$. Let $U_W \subset W|_{r_{p(x)} \bZ}$ be the sum of all subrepresentations isomorphic to $U$, and $U_W'$ the sum of all irreducible subrepresentations not isomorphic to $U$. Then the decomposition $W|_{r(p) \bZ} = U_W \oplus U_{W}'$ is canonical, and so there is a canonical projection
\[
\End_L^{\sigma}(W|_{r(p) \bZ}) \rightarrow \End_L^{\sigma^{r(p)}}(U_W).
\]
Let $\psi_{p(x)}$ denote the precomposition with the restriction map $\End_L^{\sigma}(W) \rightarrow \End_L^{\sigma^{r(p)}}(W|_{r(p) \bZ})$. 

\begin{lemma}\label{lem:psipxisom}
	The map $\psi_{p(x)} \colon \End_L^{\sigma}(W) \rightarrow \End_L^{\sigma^{r(p)}}(U_W)$ is a $K(\alpha)$-algebra isomorphism.
\end{lemma}

Here $\End_L^{\sigma}(W)$ as viewed as a $K(\alpha)$-algebra, and $\End_L^{\sigma^{r(p)}}(U_W)$ as a $K_{p(x)}(\alpha)$-algebra, both as in Lemma \ref{lem:centredesc}. The latter field $K_{p(x)}(\alpha)$ is nothing but $K(\alpha)$ by Lemma \ref{lem:dontneedpx}.

\begin{proof}
	As $\End_L^{\sigma}(W)$ is a division algebra, it suffices to show that both sides have the same $K$-dimension. We may compute using \cite[Thm.\ A]{TAYSchur} that
	\begin{align*}
	\dim_K \End_L^{\sigma}(W) &= \dim_L \End_L^{\sigma}(W|_{r \bZ}), \\
	 &= [\Gamma : \Gamma_{p(x)}] \cdot m_L^{\sigma}(p(x))^2 \cdot \dim_L \End_{L[r \bZ]}(V_{p(x), 1, r}), \\
	&= [\Gamma : \Gamma_{p(x)}] \cdot m_L^{\sigma}(p(x))^2 \cdot \deg(p),
	\end{align*}
	and
	\begin{align*}
	\dim_K \End_L^{\sigma^{r(p)}}(U_W) &= [K_{p(x)} : K] \cdot \dim_{K_{p(x)}} \End_L^{\sigma^{r(p)}}(U_W), \\
	&= [K_{p(x)} : K] \cdot m_L^{\sigma^{r(p)}}(p(x))^2 \cdot \dim_{L} \End_{L[r \bZ]}(U_W|_{r \bZ}), \\
	&= [K_{p(x)} : K] \cdot m_L^{\sigma^{r(p)}}(p(x))^2 \cdot \deg(p),
	\end{align*}
	and note that $[\Gamma : \Gamma_{p(x)}] = [K_{p(x)} : K]$ and $m_L^{\sigma}(p(x)) = m_L^{\sigma^{r(p)}}(p(x))$ by Lemma \ref{lem:reductiontopolyoverK}.
\end{proof}

We can now describe the class of the division algebra $\End_{L}^{\sigma}(W_{p(x)})$ in $\Br(K(\alpha))$.

\begin{cor}\label{cor:descofendoinbrgroup}
	Suppose that $p(x) \in \sP_L$. Then there is a $K(\alpha)$-algebra isomorphism
	\[
	(L(\alpha) / K(\alpha), \tau_{\alpha}, \alpha) \cong M_{\frac{r}{r(p) m_L^{\sigma}(p(x))}}(\End_{L}^{\sigma}(W_{p(x)})^{\opp}).
	\]
	In particular,
	\[
	[\End_{L}^{\sigma}(W_{p(x)})] = -[(L(\alpha) / K(\alpha), \tau_{\alpha}, \alpha)]
	\]
	in $\Br(K(\alpha))$.
\end{cor}

\begin{proof}
	To ease notion, set $m \coloneqq m_L^{\sigma}(p(x))$, and recall that we write $n = \deg(p)$. Using Lemma \ref{lem:mdivkmatrixform}, pick $A \in \GL_m(L(\alpha))$ with $N_m^{\tau_{\alpha}}(A) = \alpha I_m$, which defines the irreducible semilinear representation $L(\alpha)_A^{m, \tau_{\alpha}}$, for which by Lemma \ref{lem:linearbrauerequiv} there is a $K(\alpha)$-linear isomorphism 
	\[
	(L(\alpha) / K(\alpha), \tau_{\alpha}, \alpha) \cong M_{\frac{r}{r(p) m_{L(\alpha)}^{\tau_{\alpha}}(x - \alpha)}}(\End_{L(\alpha)}^{\tau_{\alpha}}(L(\alpha)_A^{m, \tau_{\alpha}})^{\opp})
	\]
	and thus, as $m_{L(\alpha)}^{\tau_{\alpha}}(x - \alpha) = m_L^{\sigma}(p(x)) = m$ from Theorem \ref{thm:schurindexpx}, it suffices to show that
	\[
	\End_{L(\alpha)}^{\tau_{\alpha}}(L(\alpha)_A^{m, \tau_{\alpha}}) \cong \End_L^{\sigma^{r(p)}}(U_W)
	\]
	as $K(\alpha)$-algebras by Lemma \ref{lem:psipxisom}. To this end, recall from Lemma \ref{lem:imageoffield1dim} that we have the injection
	\[
	M(-) \colon \colon M_m(L(\alpha)) \hookrightarrow M_{mn}(L)
	\]
	which satisfies $M(\tau_{\alpha}(X)) = \sigma^{r(p)}(M(X))$, and thus $M(-)$ restricts to an injection
	\[
	\Phi \colon \End_{L(\alpha)}^{\tau_{\alpha}}(L(\alpha)_A^{m, \tau_{\alpha}}) \rightarrow \End_L^{\sigma^{r(p)}}(L_{M(A)}^{nm, \sigma^{r(p)}}).
	\]
	We show this is an isomorphism by comparing dimensions $K$-dimensions. The left-hand side has $K(\alpha)$-dimension $m_{L(\alpha)}^{\tau_{\alpha}}(x - \alpha)^2$, hence $K$-dimension $[\Gamma : \Gamma_{p(x)}] \deg(p) m^2$, using that $[\Gamma : \Gamma_{p(x)}] = [K_{p(x)} : K]$, $\deg(p) = [K_{p(x)}(\alpha) : K_{p(x)}] = [K(\alpha) : K_{p(x)}]$, and $m = m_{L(\alpha)}^{\tau_{\alpha}}(x - \alpha)$ by Theorem \ref{thm:schurindexpx}. For the right-hand side, we note that as in the proof of Lemma \ref{lem:mdivkmatrixform},
	\[
	L_{M(A)}^{nm, \sigma^{r(p)}} \cong U_{p(x)},
	\]
	where $U_{p(x)}$ denotes the irreducible representation of $r(p) \bZ$ which restricts to $V_{p(x), 1, r}^{\oplus m}$. We have that $U_{p(x)} \cong U_{W}$, the isotypic component of $W|_{r(p) \bZ}$ considered above, and thus 
	\[
	\dim_K \End_L^{\sigma^{r(p)}}(L_{M(A)}^{nm, \sigma^{r(p)}}) = \dim_K \End_L^{\sigma^{r(p)}}(U_{W_{p(x)}}) = [\Gamma : \Gamma_{p(x)}] \cdot \deg(p) \cdot m^2
	\]
	as in the proof of Lemma \ref{lem:psipxisom}. Thus $\Phi$ is an isomorphism, and that we are done because
	\[
	\End_L^{\sigma^{r(p)}}(U_W) \cong \End_L^{\sigma^{r(p)}}(L_{M(A)}^{nm, \sigma^{r(p)}}).\qedhere
	\]
\end{proof}

\section{Representatives for the Semilinear Normal Form}\label{sect:representatives}

In this section we consider how the results of Section \ref{sect:schurindicesresults} interact with the construction and canonicity of the representatives $A_{\OO}$ for the semilinear normal form of Section \ref{sect:backgroundsemilinearreps}.

The following operation will be useful.

\begin{defn}
Let, for $m, s \geq 1$ and $B \in M_{m}(L)$, $R_s(B) \in M_{ms}(L)$ be the matrix
\[
R_s(B) \coloneqq \begin{pmatrix}
    0 & 0 &\cdots & 0 & B \\
    I_{m} & 0 & \cdots & 0 & 0 \\
    \vdots & \vdots & \ddots & \ddots & \vdots \\
    0 & 0 & \cdots & 0 & 0 \\
    0 & 0 & \cdots & I_{m} & 0
\end{pmatrix}.
\]
\end{defn} 
The key property that we will use of $R_s(B)$ is that
\[
R_s(B) \sigma(R_s(B)) \cdots \sigma^{s-1}(R_s(B)) = B \oplus \sigma(B) \oplus \cdots \oplus \sigma^{s-1}(B).
\]

\subsection{Construction of Representatives $A_{\OO}$}

Suppose that $\OO \in \Gamma \backslash (\sP_L \times \bZ_{\geq 1})$. In this section we consider the question of how to construct an explicit $A_{\OO} \in \GL_{N}(L)$ with $L^{N, \sigma}_{A_{\OO}} \cong W_{\OO}$, thus providing an explicit model for the indecomposable representation $W_{\OO}$ (recall that $W_{\OO}$ denotes the indecomposable representation of $\bZ$ over $L$ corresponding to the orbit $\OO$).

We make a series of reductions, which correspond exactly to the reductions that were made to compute the semilinear Schur indices outlined at the start of Section \ref{sect:schurindicesresults}.

First, let us pick a representative $p(x) \in \sP_L$ so that $\OO = [(p(x), t)]$, let $m \coloneqq m_L^{\sigma}(p(x)^t)$, and note that the $L$-dimension $W_{\OO}$ is $N \coloneqq mnr(p)t$, where $n = \deg(p)$ and $r(p) = [K_{p(x)} : K]$. As in the proof of Lemma \ref{lem:mdivkmatrixform}, the condition that $L^{N, \sigma}_{A_{\OO}} \cong W_{\OO}$ is satisfied whenever
\[
N^{\sigma}_{N}(A_{\OO}) = \bigoplus_{\gamma \in \Gamma / \Gamma_{p(x)}} C(\gamma * p^t)^{\oplus m},
\]
and so this is our goal.

\subsubsection*{First Reduction} Suppose that we have $B(p^t) \in \GL_{mnt}(L)$ with 
\[
N^{\sigma^{r(p)}}_{mnt}(B(p^t)) = C(p^t)^{\oplus m}.
\] 
Then setting
\[
A_{\OO} \coloneqq R_{r(p)}(B(p^t)),
\]
we have that
\begin{align*}
N^{\sigma}_{N}(A_{\OO}) &= N^{\sigma^{r(p)}}_{N}(B(p^t) \oplus \sigma(B(p^t)) \oplus \sigma^{r(p) - 1}(B(p^t))), \\
&= N^{\sigma^{r(p)}}_{mnt}(B(p^t)) \oplus N^{\sigma^{r(p)}}_{mnt}(\sigma(B(p^t))) \oplus \cdots \oplus N^{\sigma^{r(p)}}_{mnt}(\sigma^{r(p) - 1}(B(p^t))), \\
&=  C(p^t)^{\oplus m} \oplus C(\sigma * p^t)^{\oplus m} \oplus \cdots \oplus C(\sigma^{r(p) - 1} * p^t)^{\oplus m},
\end{align*}
as required, and thus one is reduced to finding such a $B(p^t) \in \GL_{mnt}(L)$.

\subsubsection*{Second Reduction} Let us first focus on finding $B(p^1)$. Fix a root $\alpha \in \overline{L}$ of $p(x)$, and suppose that $B_{\alpha}(p) \in \GL_{m}(L(\alpha))$ with
\[
N_{m}^{\tau_{\alpha}}(B_{\alpha}(p)) = \alpha I_m.
\]
Using the map $M(-) \colon \GL_m(L(\alpha)) \rightarrow \GL_{mn}(L)$, 
\[
N_{mn}^{\sigma^{r(p)}}(M(B_{\alpha}(p))) = M(N_{m}^{\tau_{\alpha}}(B_{\alpha}(p))) = M(\alpha I_m) = C(p)^{\oplus m},
\]
and so one may take
\[
B(p^1) \coloneqq M(B_{\alpha}(p)).
\]

\subsubsection*{Third Reduction} Suppose that we have found such a $B_{\alpha}(p)$, and hence $B(p^1)$. We can then construct $B(p^t)$ for all $t \geq 1$ as follows, following the proof of Theorem \ref{thm:mpequalsmpt}.

Write $B_{\alpha}(p) = (b_{ij})$, where $b_{ij} = q_{ij}(\alpha) \in L(\alpha)$ is written with respect to the basis $1, \alpha, ... \alpha^{n-1}$ of $L(\alpha)$, so each $q_{ij} \in L[y]/p(y)$.  We have for $t \geq 1$, $L$-algebras $R_t \coloneqq L[y] / p(y)^t$. Fix $u \in L / K_{p(x)}$ with $\Tr(u) = 1$, possible as $L / K$ is separable. For $t \geq 1$, we inductively construct $A_t \in \GL_m(R_t)$ as follows. When $t = 1$, we take $A_1 \coloneqq (q_{ij}(y)) \in \GL_m(R_1)$. Then for $t \geq 1$, given $A_t \in \GL_m(R_t)$,
\begin{enumerate}
	\item Take any lift $A_{t+1,0} \in M_m(R_{t+1})$, which will automatically be invertible,
	\item Set
	\[
		E_{t+1} \coloneqq A_{t+1,0} \cdot \sigma^{r(p)}(A_{t+1, 0}) \cdot \cdots  \cdot \sigma^{(\frac{r}{r(p)}-1)r(p)}(A_{t+1, 0}) - y I_m \in M_m(R_{t+1}),
	\]
	\item Define
	\[
	A_{t+1} \coloneqq (1 - u y^{-1}E_{t+1})A_{t+1, 0} \in \GL_m(R_{t+1}).
	\]
\end{enumerate}

Then, if $A_{t+1} = (f_{ij}) \in \GL_m(R_t)$, we may take $B(p^t)$ to be the block matrix
\[
B(p^t) \coloneqq (f_{ij}(C(p^t))) \in \GL_{mnt}(L).
\]

\subsubsection*{Summary} We have thus reduced the construction of $A_{\OO}$ to the construction of $B_{\alpha}(p) \in \GL_{m}(L(\alpha))$ with $N_{m}^{\tau_{\alpha}}(B_{\alpha}(p)) = \alpha I_m$. Let us know consider some cases where this is solvable.

\subsection{Example: Trivial Semilinear Schur Indices}

When $m = 1$, we can take $B_{\alpha}(p) = (\beta)$, where $\beta \in L(\alpha)$ is any element with $N_{L(\alpha) / K(\alpha)}(\beta) = \alpha$.

\subsection{Example: Maximal Semilinear Schur Indices} We know from \cite[Prop.\ 9.3(6)]{TAYSchur} that $m \mid [L(\alpha) : K(\alpha)]$. Suppose that $m = [L(\alpha) : K(\alpha)]$. Then one can take $B_{\alpha}(p) = R_m(\alpha)$.

\subsection{Example: Prime Degree Extensions} Suppose that $[L:K]$ is a prime. We know that $m \mid [L:K]$ from \cite[Prop.\ 9.3(6)]{TAYSchur}, and so if $m \neq 1$ (which is dealt with above) then $m = [L:K]$. In this case, as $m = m_{L(\alpha)}^{\tau_{\alpha}}(x - \alpha)$, we similarly have that $m \mid [L(\alpha) : K(\alpha)]$, and thus as $[L(\alpha) : K(\alpha)] = [L : K_{p(x)}] \mid [L:K]$, $[L : K_{p(x)}] = [L:K]$ and so $K = K_{p(x)}$: i.e.\ $p(x)$ is valued over $K$, and the first reduction part of the construction of $A_{\OO}$ can be avoided. More importantly, when $m = [L:K]$, we can explicitly take $B_{\alpha}(p) = R_m(\alpha)$ as above.

In summary: when $[L:K]$ is prime we can always construct the $B_{\alpha}(p)$ provided we can solve $N_{L(\alpha) / K(\alpha)}(\beta) = \alpha$ for $\beta \in L(\alpha)$. For example this is possible if $L$ is a number field \cite{Denis}.

\begin{eg}
When $L / K$ is $\bC / \bR$ we can be very explicit and give $A_{\OO}$ directly avoiding the above reductions as follows. Let us identify $\Gamma \backslash (\sP_{\bC} \times \bZ_{\geq 1})$ with $\bH_{\geq 0}^{\times} \times \bZ_{\geq 1}$, where $\bH_{\geq 0}^{\times} \coloneqq \{z \in \bC^{\times} \mid \text{im}(z) \geq 0\}$. Recall from Example \ref{eg:CoverR} that for $z \in \bC^{\times}$, $m_{\bC}^{\sigma}(x-z) = 2$ if $z \in \bR_{<0}$ and $m_{\bC}^{\sigma}(x-z) = 1$ otherwise. Therefore, for $\OO_{z,t}$ for the orbit corresponding to any $(z,t) \in \bH_{\geq 0}^{\times} \times \bZ_{\geq 1}$, as described in Section \ref{sect:backgroundsemilinearreps} a model for $W_{\OO_{z,t}}$ is given by finding a matrix $A_{z,t}$ such that the norm
\[
A_{z,t} \overline{A_{z,t}} \sim \twopartdef{C((x-z)^t)}{z \in \bR_{>0}}{C((x-z)^t) \oplus C((x-\overline{z})^t)}{z \not\in \bR_{>0}}
\]
If $z \in \bR_{> 0}$, we take $A_{z,t} \coloneqq J_t(\sqrt{z})$, and for $z \not\in \bR_{> 0}$
\[
A_{z,t} \coloneqq \begin{pmatrix}
	0 & J_t(z) \\ 1 & 0
\end{pmatrix}.
\]
We can check that when $z \in \bR_{> 0}$,
\[
A_{z,t} \overline{A_{z,t}} = J_t(\sqrt{z})^2 \sim J_t(z) \sim C((x - z)^t),
\]
and when $z \not\in \bR_{> 0}$,
\[
A_{z,t} \overline{A_{z,t}} = J_t(z) \oplus J_t(\overline{z}) \sim C((x-z)^t) \oplus C((x-\overline{z})^t
\]

Now Corollary \ref{cor:SNF} restated in this context is as follows. Suppose $n \geq 1$ and $A \in \GL_n(\bC)$. Then
there are $z_1, ... , z_a \in \bH_{\geq 0}^{\times} \setminus \bR_{>0}$, $y_1, ... , y_b \in \bR_{>0}$, and $t_1, ... , t_a, s_1, ... s_b \geq 1$ such that
\[
A \sim_{\sigma} \begin{pmatrix}
	0 & J_{t_1}(z_1) \\ 1 & 0
\end{pmatrix} \oplus \cdots \oplus \begin{pmatrix}
	0 & J_{t_a}(z_a) \\ 1 & 0
\end{pmatrix} \oplus J_{s_1}(\sqrt{y_1}) \oplus \cdots \oplus J_{s_b}(\sqrt{y_b}),
\]
and this description is unique up to permutation of the blocks. Furthermore, one can read off the $(z_i,t_i)$ and $(y_j,s_j)$ from the Jordan normal form of $A \overline{A}$, which will be
\[
A\overline{A} \sim J_{t_1}(z_1) \oplus \cdots \oplus J_{t_a}(z_a) \oplus J_{t_1}(\overline{z_1}) \oplus \cdots \oplus J_{t_a}(\overline{z_a}) \oplus J_{s_1}(y_1) \oplus \cdots \oplus J_{s_{b}}(y_b).
\]
This recovers \cite[Thm.\ 3.1]{HongHorn} in the case when $A$ is invertible.
\end{eg}

\subsection{Unramified Extensions of Local Fields}

In Section \ref{sect:localfields} below we give an explicit construction of $B_{\alpha}(p)$ when $L / K$ is an unramified extension of local fields.

\section{Description of Indecomposables and Base Change}\label{sect:basechange}

\subsection{Description of Indecomposables}\label{sect:altdescofirred}

In this section we give an alternative description of the set $\Ind_L^{\sigma}(\bZ)$, which until now we have viewed as $\Gal(L/K) \backslash (\sP_L \times \bZ_{\geq 1}) = (\Gal(L/K) \backslash \sP_L) \times \bZ_{\geq 1}$.

The following is a standard fact, for which we were unable to find a reference.

\begin{lemma}
	Suppose that $K$ is a field and $F_1, F_2$ are (potentially infinite) Galois extensions of $K$ with $F_1 \subset F_2$. Then there is a bijection
	\[
	\Gal(F_1/K) \backslash \sP_{F_1} \leftrightarrow \Gal(F_2/K) \backslash \sP_{F_2},
	\]
	where any representative $q_1 \in \sP_{F_1}$ maps to any irreducible factor of $q_1$ over $F_2$, and any representative $q_2 \in \sP_{F_2}$ maps to the product over all $\Gal(F_2 / F_1)$-conjugates, which lies in $\sP_{F_1}$.
\end{lemma}

\begin{proof}
	Let $G_1 \coloneqq \Gal(F_1 / K)$, $G_2 \coloneqq \Gal(F_2 / K)$, let $H \coloneqq \Gal(F_2 / F_1)$, and for $q_2 \in \sP_{F_2}$ set
	\[
	\Pi(q_2) \coloneqq \prod_{q \in H \cdot q_2} q.
	\]
	This product is finite, as the orbit $H \cdot q_2$ is finite, because $q_2$ has coefficients in a finite extension of $F_1$. Furthermore, $\Pi(q_2)$ has coefficients in $F_1$ because $\Pi(q_2)$ is stable under $H$ and $F_2^H = F_1$. Given a monic non-constant divisor $f \mid \Pi(q_2)$ over $F_1$, let $\OO_f \subset H \cdot q_2$ denote the irreducible factors of $f$. As $f$ has coefficients in $F_1$, $\OO_f$ is stable under the action of $H$, and therefore $\OO_f = H \cdot q_2$, hence $f = \Pi(q_2)$. Thus $\Pi(q_2)$ is monic and irreducible, hence $\Pi(-)$ is a well-defined map $\sP_{F_2} \rightarrow \sP_{F_1}$.
	
	To see that this induces a map $\Phi \colon G_2 \backslash \sP_{F_2} \rightarrow G_1 \backslash \sP_{F_1}$, note that for any $\sigma \in G_2$,
	\[
	\Pi(\sigma(q_2)) = \prod_{q \in H \cdot \sigma(q_2)} q = \prod_{q \in \sigma \cdot (\sigma^{-1} H \sigma) \cdot q_2} q
	= \prod_{q \in \sigma H \cdot q_2} q = \prod_{q \in H \cdot q_2} \sigma(q) = \sigma'(\Pi(q_2)),
	\]
	where $\sigma' \in G_1$ is the image of $\sigma$, and we have used that $H$ is normal in $G_2$.

	In the other direction, given $q_1 \in \sP_{F_1}$, note that for any monic irreducible factor $p$ of $q_1$ over $F_2$, $q_1 = \Pi(p)$. Indeed, for any $\sigma \in H$, $\sigma(p) \mid q_1$, and so $\Pi(p) \mid q_1$, hence by irreducibility $q_1 = \Pi(p)$. In particular, the map $\sP_{F_1} \rightarrow G_2 \backslash \sP_{F_2}$ is independent of the choice of irreducible factor. Using this, we can show that this map descends to a map $\Psi \colon G_1 \backslash \sP_{F_1} \rightarrow G_2 \backslash \sP_{F_2}$. Indeed, given $\sigma_1 \in G_1$, take a lift $\sigma_2 \in G_2$: for $q_1 = p_1 \cdots p_m$ over $F_2$, $\sigma_1(q_1) = \sigma_2(q_1) = \sigma_2(p_1) \cdots \sigma_2(p_m)$, and so, as we may freely choose the irreducible factor, $\sigma_1(q_1)$ maps to $\sigma_2(p_1)$.

	Finally we note that both compositions of $\Phi$ and $\Psi$ are the identity. Indeed, for $q_2 \in \sP_{F_2}$, $\Phi(q_2) = \Pi(q_2)$, of which $q_2$ is an irreducible factor over $F_2$, hence $\Psi(\Phi(q_2)) = q_2$. For $q_1 \in \sP_{F_1}$, we have that $q_1 = \prod_{q \in H \cdot p} q$ for any irreducible factor $p$ of $q_1$ over $F_2$, and so $q_1 = \Phi(p) = \Phi(\Psi(q_1))$.
\end{proof}

We are interested because it allows us to give three different points of view on $\Ind_L^{\sigma}(\bZ)$. Firstly, we may take $F_1 = K$ and $F_2 = L$ to obtain a canonical identification
\[
\Ind_L^{\sigma}(\bZ) \leftrightarrow (\Gal(L/K) \backslash \sP_L) \times \bZ_{\geq 1} \leftrightarrow \sP_K \times \bZ_{\geq 1}.
\]

Explicitly, for $(q,t) \in \sP_K \times \bZ_{\geq 1}$, the corresponding indecomposable semilinear representation is determined by the property that the restriction to $r \bZ$ is a direct sum of copies of the $r \bZ$-representation $L[x] / q(x)^t$, where the generator $r \in r\bZ$ acts by multiplication by $x$.

Secondly, if $K$ is perfect, then $\overline{K} / K$ is Galois, and we may take $F_1 = L$ and $F_2 = \overline{K}$. Therefore, writing $G_K \coloneqq \Gal(\overline{K} / K)$, we have a canonical identification 
\[
\Ind_L^{\sigma}(\bZ) \leftrightarrow (\Gal(L/K) \backslash \sP_L) \times \bZ_{\geq 1} \leftrightarrow \Gal_K \! \backslash \overline{K}^{\times} \times \bZ_{\geq 1}.
\] 
as $\sP_{\overline{K}}$ is canonically identified with $\overline{K}^{\times}$. Explicitly, a given $(\alpha,t) \in \overline{K}^{\times} \times \bZ_{\geq 1}$, $\alpha$ corresponds to $W_{m_{\alpha, L}, t} \in \Ind_{L}^{\sigma}(\bZ)$, where $m_{\alpha, L} \in \sP_L$ denotes the minimum polynomial of $\alpha$ over $L$. If $m_{\alpha, K} \in \sP_K$ denotes the minimum polynomial of $\alpha$ over $K$, then $W_{m_{\alpha, L}, t}|_{r \bZ} \cong (L[x] / m_{\alpha, K}(x)^t)^{m_L^{\sigma}(m_{\alpha, L}(x))}$, where $r$ acts by multiplication by $x$.
We use this second description to describe base change below.

\subsection{General Base Change}

Suppose that $G$ is a group, $L_1 \subset L_2$ are fields, and $\sigma_1 \colon G \rightarrow \Aut(L_1)$, $\sigma_2 \colon G \rightarrow \Aut(L_2)$ are group homomorphisms such that $\sigma_2(g)|_{L_1} = \sigma_1(g)$ for all $g \in G$.

In this context, there is a base change functor
\[
L_2 \otimes_{L_1} (-) \colon \Rep_{L_1}^{\sigma_1}(G) \rightarrow \Rep_{L_2}^{\sigma_2}(G)
\]
where for $W \in \Rep_{L_1}^{\sigma_1}(G)$, $G$ acts on $L_2 \otimes_{L_1} W$ diagonally.

\begin{eg}
If $L_2 \coloneqq L$, and $L_1 \coloneqq K = L^G$, then $L \otimes_K (-)$ is exactly the base change functor considered in \cite[\S 4]{TAYSchur}. In particular, $L_2 \otimes_{L_1} (-)$ is neither full nor faithful in general.
\end{eg}

\subsection{Base Change for Representations of $\bZ$}

Now take $G = \bZ$. In this case, we suppose $L_1 \subset L_2$ are fields, $\sigma_1 \in \Aut(L_1)$ and $\sigma_2 \in \Aut(L_2)$, and $\sigma_2|_{L_1} = \sigma_1$. Suppose additionally that $\sigma_1$ has order $r_1$ and $\sigma_2$ has order $r_2$. Note that $\sigma_1^{r_2} = 1$, and thus $r_1 \mid r_2$.

In this setting it is natural to ask what the image is of each indecomposable modules $W_{p(x),t}$ under $L_2 \otimes_{L_1} (-)$, given as a sum of indecomposable modules.

For simplicity we assume that $L_1 / L_2$ is Galois, $L_1^{\sigma_1} = L_2^{\sigma_2}$ and that this $K \coloneqq L_1^{\sigma_1}$ is perfect. This is enough generality to cover the case of isocrystals which we consider in Section \ref{sect:basechangeisoc}. The method described here extends to the case when these assumptions are removed.

Because $K$ is perfect, we can associate to any $\alpha \in \overline{K}^{\times}$ and $t \geq 1$ an irreducible representation $W_{\alpha, t, 1} \in \Ind_{L_1}^{\sigma_1}(\bZ)$, which depends only on the $\Gal_K$-orbit of $\alpha$ and induces an identification
\[
\left( \Gal_K \! \backslash \overline{K}^{\times} \right) \times \bZ_{\geq 1} \xrightarrow{\sim} \Ind_{L_1}^{\sigma_1}(\bZ), \qquad (\alpha, t) \mapsto W_{\alpha, t, 1}.
\]
This follows from Section \ref{sect:altdescofirred}. Explicitly, $W_{\alpha, t, 1}$ is the representation $W_{p(x), t} \in \Ind_{L_1}^{\sigma_1}(\bZ)$ in our usual notation, where $p(x) \coloneqq m_{\alpha, L}$ is the minimal polynomial of $\alpha$ over $L$. If $m_{\alpha, K}$ denotes the minimal polynomial of $\alpha$ over $K$, then $W_{\alpha, t, 1}$ restricts to a direct sum of copies of the $r_1 \bZ$-representation $L[x] / m_{\alpha, K}(x)^t$, where $r_1$ acts by multiplication by $x$. We similarly have an identification
\[
\left( \Gal_K \! \backslash \overline{K}^{\times} \right) \times \bZ_{\geq 1} \xrightarrow{\sim} \Ind_{L_2}^{\sigma_2}(\bZ), \qquad (\alpha, t) \mapsto W_{\alpha, t, 2}.
\]

In what follows, we write $s \coloneqq \frac{r_2}{r_1}$, and $m_{L_1}^{\sigma_1}(\alpha, t) \coloneqq m_{L_1}^{\sigma_1}(m_{\alpha, L_1}^t)$, $m_{L_2}^{\sigma_2}(\alpha^s, t) \coloneqq m_{L_2}^{\sigma_2}(m_{\alpha^s, L_2}^t)$ for the semilinear Schur indices corresponding to $(\alpha, t)$ in $\Ind_{L_1}^{\sigma_1}(\bZ)$ and to $(\alpha^s, t)$ in $\Ind_{L_2}^{\sigma_2}(\bZ)$.

\begin{prop}\label{prop:explicitdescofbasechange}
	Suppose that $\alpha \in \overline{K}^{\times}$, and $t \geq 1$. Then if $t = 1$ or $\charfield(K) \nmid s$,
	\[
		L_2 \otimes_{L_1} W_{\alpha, t, 1} \cong W_{\alpha^s, t, 2}^{\oplus N}, \: \text{ where } \:  N \coloneqq \frac{[K(\alpha) : K(\alpha^s)] \cdot m_{L_1}^{\sigma_1}(\alpha, t)}{m_{L_2}^{\sigma_2}(\alpha^s, t)}.
	\]

	In particular, $m_{L_2}^{\sigma_2}(\alpha^s, t) \mid [K(\alpha) : K(\alpha^s)] \cdot m_{L_1}^{\sigma_1}(\alpha, t)$.
\end{prop}

\begin{proof}
	Let $y \coloneqq x^s$, and identify representations of $r_1 \bZ$ with $L_1[x]$-modules, so that restriction of $L_1$-linear representations of $r_1 \bZ$ to $r_2 \bZ$ corresponds to restricting the corresponding $L_1[x]$-module to $L_1[y]$. As 
	\[
	W_{\alpha^s, t, 2}|_{r_2 \bZ} \cong (L_2[y] / m_{\alpha^s, K}(y)^t)^{\oplus m_{L_2}^{\sigma_2}(\alpha^s, t)},
	\]
	it suffices to show that
	\[
	(L_2 \otimes_{L_1} W_{\alpha, t, 1})|_{r_2 \bZ} \cong (L_2[y] / m_{\alpha^s, K}(y)^t)^{\oplus [K(\alpha) : K(\alpha^s)] \cdot m_{L_1}^{\sigma_1}(\alpha, t)}
	\]
	by Theorem \ref{thm:isomonrestrictions}. The left-hand side is $(L_2 \otimes_{L_1} W_{\alpha, t, 1}|_{r_1 \bZ})|_{r_2 \bZ}$, and
	\[
	W_{\alpha, t, 1}|_{r_1 \bZ} \cong (L_1[x] / m_{\alpha, K}(x)^t)^{\oplus m_{L_1}^{\sigma_1}(\alpha, t)},
	\]
	and thus it suffices to show that the $L_2[x]$-module
	\[
	L_2 \otimes_{L_1} L_1[x] / m_{\alpha, K}(x)^t = L_2[x] / m_{\alpha, K}(x)^t
	\]
	restricts as an $L_2[y]$-module to
	\[
	(L_2[y] / m_{\alpha^s, K}(y)^t)^{\oplus [K(\alpha) : K(\alpha^s)]}.
	\]
	For this, as both polynomials are defined over $K$, it suffices to show that
	\[
	 (K[x] / m_{\alpha, K}(x)^t)|_{K[y]} \cong (K[y] / m_{\alpha^s, K}(y)^t)^{\oplus [K(\alpha) : K(\alpha^s)]}.
	\]

	When $t = 1$, this is clear, as $K[x] / m_{\alpha, K}(x) = K(\alpha)$ and $K[y] / m_{\alpha^s, K}(y) = K(\alpha^s)$, and this isomorphism is simply writing $K(\alpha)$ as a $K(\alpha^s)$-vector space.

	For general $t \geq 1$, with $\charfield(L) \nmid s$, let $p(x) \coloneqq m_{\alpha, K}(x) \in K[x]$ and $q(y) \coloneqq m_{\alpha^s, K}(y) \in K[y]$, monic and irreducible. As $q(\alpha^s) = 0$, $q(x^s)$ has $\alpha$ as a root, hence $p(x) \mid q(x^s)$. We claim this is a simple root, which implies that $q(x^s) = p(x) h(x)$, where $h \in K[x]$ is coprime to $p(x)$. Indeed, if we write $q(y) = (y- \beta_1) \cdots (y - \beta_m)$ for $\beta_i \in \overline{K}$, which are distinct as $K$ is perfect, fixing $\beta_1 = \alpha^s$, $q(x^s) = (x^s- \alpha^s) \cdots (x^s - \beta_m)$. The $\beta_i$ are distinct, so $(\alpha^s - \beta_i) \neq 0$ for $i \geq 2$, and $x^s - \alpha^s$ has $\alpha$ as a simple root because the formal derivative $s x^{s-1}$ is non-zero at $\alpha$, as $s \neq 0$ in $K$ and $\alpha \neq 0$.  

	Now $M \coloneqq (K[x] / p(x)^t)|_{K[y]}$ is killed by $q(y)^t$, hence 
	\[
		M \cong \bigoplus_{i = 1}^d K[y] / q(y)^{e_i}
	\]
	for some $d, e_i \geq 1$ with each $1 \leq e_i \leq t$. Now we know that $q(x^s)M = p(x)M$ by the above, hence 
	\[
		d = \dim_{K[y] / q(y)} M / q(y) M = \dim_{K[y] / q(y)} M / p(x) M = [K(\alpha) : K(\alpha^s)].
	\]
	We also know that
	\[
		\sum_{i = 1}^d e_i \deg(q) = \dim_K M = t \deg(p) = t [K(\alpha) : K] = d t\deg(q),
	\]
	and each $e_i \leq t$, which forces each $e_i = t$ and the statement follows. 
\end{proof}

\section{Applications to the Category of Isocrystals Over a Finite Field}\label{sect:localfields}

In this section we focus on the case that $L$ is a local field. We then specialise to the case when $L / K$ is unramified, and then again to the case of isocrystals over a finite field.

\subsection{Non-Archimedean Local Fields}

Suppose that $K$ is a non-archimedean local field, with residue field of order $q$. For an unramified extension $L / K$, we write $\Frob_K \in \Gal(L / K)$ for the automorphism of $L$ which induces the automorphism $x \mapsto x^q$ on the residue field.

We write $\inv_K \colon \Br(K) \xrightarrow{\sim} \bQ / \bZ$ for the Haase invariant, as defined in \cite[\S IV.4]{milneCFT}. This is characterised by taking the values $\inv_K([(K_n / K, \Frob_K, \pi)]) = 1/n$ for any $n \geq 1$ and uniformiser $\pi$ of $K$. Let $\phi_{K} \colon K^{\times} \rightarrow \Gal(K^{\ab} / K)$ denote the local Artin map, and for any finite abelian extension $L / K$ write $\phi_{L / K}$ for the composition of $\phi_K$ with $ \Gal(K^{\ab} / K) \rightarrow \Gal(L/K)$. Here $\phi_K$ is normalised so that $\phi_{L/K}(\pi) = \Frob_K$ whenever $L / K$ is unramified, as in \cite[\S III.3]{milneCFT}.

We have the following immediate description of the Hasse invariant of $\End_{L \rtimes_{\sigma} \bZ}(W_{\lambda})$.

\begin{prop}\label{prop:brauereltdesc}
Suppose $K$ is a non-archimedean local field, and $\lambda \in K^{\times}$. Let $\chi \colon \Gal(L/K) \rightarrow \bQ / \bZ$ denote the character with $\chi(\sigma) = 1/r$. Then
\[
\inv_K(\End_{L \rtimes_{\sigma} \bZ}(W_{\lambda})) = -\chi(\phi_{L/K}(\lambda)).
\]
\end{prop}

\begin{proof}
	Writing $\psi_L \colon H^2(L/K) \xrightarrow{\sim} \Br(L/K)$ for the isomorphism of \cite[Thm.\ IV. 3.14]{milneCFT}, we have that $\chi(\phi_{L/K}(\lambda)) = \inv_K(\psi_L(a \cup \delta \chi))$ by \cite[Prop.\ III.3.6 ]{milneCFT} and \cite[Rem.\ IV.4.4(a)]{milneCFT}. By definition, the $\bZ$-valued $2$-cocycle $d\chi$ on $\Gal(L/K)$ has values
	\[
	d\chi(\sigma^i, \sigma^j) = \twopartdef{0}{i+j < r}{1}{i+j \geq r}
	\]
	and thus $c \coloneqq a \cup \delta \chi$ has values
	\[
	c(\sigma^i, \sigma^j) = \twopartdef{1}{i+j < r}{\lambda}{i+j \geq r}
	\]
	as in the discussion before \cite[Prop.\ III.3.6 ]{milneCFT}, hence $\psi_L(a \cup \delta \chi)$ is nothing but $[(L/K, \sigma, \lambda)]$. The claim then follows from Lemma \ref{lem:linearbrauerequiv}, that $[\End_{L \rtimes_{\sigma} \bZ}(W_{\lambda})] = -[(L/K, \sigma, \lambda)]$. 
\end{proof}

Pairing this with Proposition \ref{prop:brauereltdesc}, we also can describe the Hasse invariant of $\End_{L \rtimes_{\sigma} \bZ}(W_{p(x),1})$, which we may canonically view as a central division algebra over $K(\alpha)$ in light of Section \ref{sect:endomorphimrings}.

\begin{cor}
Suppose that $K$ is a non-archimedean local field and $p(x) \in \sP_L$, and let $\alpha \in \overline{L}$ be a root of $p(x)$. Then for the character $\chi \colon \Gal(L(\alpha) / K(\alpha)) \rightarrow \bQ / \bZ$ with $\chi(\tau_{\alpha}) = \frac{r(p)}{r}$, 
\[
\inv_{K(\alpha)}(\End_{L \rtimes_{\sigma} \bZ}(W_{p(x),1})) = -\chi(\phi_{L(\alpha) / K(\alpha)}(\alpha)).
\]
\end{cor}

\subsection{Unramified Extensions of Local Fields}\label{sect:unramexts}

From Proposition \ref{prop:brauereltdesc}, we immediately obtain the following simple formulas for semilinear Schur indices when $L / K$ is unramified.

\begin{cor}\label{cor:unramSLSchur}
	If $L / K$ is unramified and $\sigma = \Frob_K$, then for any $\lambda \in K^{\times}$, 
	\[
	\inv_K(\End_{L \rtimes_{\sigma} \bZ}(W_{\lambda})) = -\frac{\nu_K(\lambda)}{r}.
	\]
	In particular,
	\[
	m_L^{\sigma}(x - \lambda) = \frac{r}{\gcd(r, \nu_K(\lambda))}.
	\]
\end{cor}

\begin{proof}
	When $L / K$ is unramified, $\phi_{L/K}(\lambda) = \Frob_K^{\nu_K(\lambda)}$ by \cite[Lem.\ III.3.7]{milneCFT}, and thus the first claim follows directly from Proposition \ref{prop:brauereltdesc}. The second follows from the fact that $m_L^{\sigma}(x - \lambda)$ is the order of $\End_{L \rtimes_{\sigma} \bZ}(W_{\lambda})$ in $\Br(K)$ by Corollary \ref{cor:cyclicindex} and \cite[Rem.\ IV.4.4(b)]{milneCFT}.
\end{proof}

\begin{remark}
If one is only interested in description of $m_L^{\sigma}(x - \lambda)$, this can be prove more simply, using the norm and valuation maps.

We write $L^\times  = \langle \pi \rangle \times \OO_L^\times$ and $K^\times  = \langle \pi \rangle \times \OO_K^\times$, where $\pi$ is a uniformiser of $K$. Because the restricted norm map $N_{L/K}: \OO_L^\times  \rightarrow \OO_K^\times$ is surjective \cite[Prop.\ III.1.2]{milneCFT}, the image $N_{L/K}(L^\times)$ is
\begin{align*}
    N_{L/K}(L^\times) & = N_{L/K}(\langle \pi \rangle) \times N_{L/K}(\OO_L^\times) \\
    & = \langle \pi^r \rangle \times \OO_K^\times,
\end{align*}
and so the valuation map $\nu_K: K^\times \rightarrow \bZ$ induces an isomorphism
\[
\psi \colon K^{\times} / N_{L/K}(L^{\times}) \xrightarrow{\sim}  \bZ / r \bZ, \ \ \ \ \ \lambda \mapsto \nu_K(\lambda).
\]
By Theorem \ref{thm:mainthm}, $m_L^{\sigma}(x - \lambda)$ is the order of $\lambda $ in  $K^{\times} / N_{L/K}(L^{\times})$, and so the isomorphism $\psi$ gives that $m_L^{\sigma}(x - \lambda)$ is the order of the valuation $\nu(\lambda)$ in $\bZ/r\bZ$, which is $\frac{r}{\gcd(r, \nu_K(\lambda))}$.
\end{remark}

\begin{cor}\label{cor:unramSLSchur2}
    If $L / K$ is unramified and $\sigma = \Frob_K$, then for any $p(x) \in \sP_L$ and $t \geq 1$,
    \[
        m_L^\sigma(p(x)^t) = \frac{s}{\gcd(s, \nu_{K(\alpha)}(\alpha))}
    \]
    where $\alpha$ is any root of $p(x)$ in $\overline{L}$, and $s \coloneqq [L : K_{p(x)}]$.
\end{cor}

\begin{proof}
	For any $p(x) \in \sP_L$, writing $r(p) = [K_{p(x)}:K]$, from Theorem \ref{thm:mainthm} we have that
	\[
    m_L^{\sigma}(p(x)^t) = m_{L(\alpha)}^{\tau_\alpha}(x - \alpha).
    \]
	As $L(\alpha)$ is the compositum of $L$ and $K(\alpha)$, the extension $L(\alpha)/K(\alpha)$ is unramified by \cite[Prop.\ 7.2]{Neu}. Thus from Corollary \ref{cor:unramSLSchur},
	\[
	m_L^{\sigma}(p(x)^t) = m_{L(\alpha)}^{\tau_\alpha}(x - \alpha) = \frac{s}{\gcd(s, \nu_{K(\alpha)}(\alpha))},
	\]
	where now $s$ is the degree of the unramified extension $L(\alpha)/K(\alpha)$.
\end{proof}

In the unramified case we may also describe the representatives $A_{\OO}$ in the semilinear normal form of Section \ref{sect:backgroundsemilinearreps} explicitly. As described in Section \ref{sect:representatives}, to do this it is sufficient to find, for $p(x) \in \sP_L$ and $\alpha \in \overline{L}$ a root of $p(x)$, a matrix $B_{\alpha}(p) \in \GL_m(L(\alpha))$, where $m \coloneqq m_L^{\sigma}(p(x))$, with 
\[
N^{\tau_{\alpha}}_m(B_{\alpha}(p)) = \alpha I_m.
\]
To do this: 
\begin{enumerate}
	\item Write $\alpha = \pi^{\nu} u$, where $\pi$ is a uniformiser of $K(\alpha)$, $\nu \in \bZ$ and $u \in \OO_{K(\alpha)}^{\times}$,
	\item Let $w \in \OO_{L(\alpha)}^{\times}$ with $N_{L(\alpha) / K(\alpha)}(w) = u$, possible because $L(\alpha) / K(\alpha)$ is unramified.
	\item Then for $s \coloneqq [L : K_{p(x)}]$, set $g_{\alpha} \coloneqq \gcd(s,\nu)$ and
\[
B_{\alpha}(p) \coloneqq w R_m(\pi^{\frac{\nu}{g_{\alpha}}}) \in \GL_{m_{\alpha}}(L(\alpha)),
\]
where $R_m(-)$ is a defined in Section \ref{sect:representatives}.
\end{enumerate}

Then we may check, using that $m = s / g_{\alpha}$, that 
\[
N^{\tau_{\alpha}}_m(B_{\alpha}(p)) = N_{L/K}(w) R_m(\pi^{\frac{\nu}{g_{\alpha}}})^{s} = u (\pi^{\frac{\nu}{g_{\alpha}}})^{\frac{s}{m}} I_m = \pi^{\nu} u I_m = \alpha I_m.
\]

\subsection{Isocrystals}\label{sect:isocrystals}

For a perfect field $k$, recall that we may consider the category $\Isoc(k)$ as defined in the introduction. As described there, when $k = \bF_q$ is a finite field for $q = p^r$ and $p$ a prime, there is a canonical identification 
\[
\Rep_{\bQ_q}^{\varphi_q}(\bZ) \xrightarrow{\sim} \Isoc(\bF_q),
\]
where $\bQ_q$ denotes the unramified extension of $\bQ_p$ of degree $r$, and $\varphi_q$ is a lift of the Frobenius $[x \mapsto x^p] \in \Gal(\bF_q / \bF_p)$. In particular, as $\bQ_q / \bQ_p$ is unramified, the results above give use an explicit description of semilinear Schur indices in this case.

The fours facts about $\Isoc(\bF_q)$ from the introduction follows from the general facts about $\Rep_{\bQ_q}^{\varphi_q}(\bZ)$ of Section \ref{sect:backgroundsemilinearreps}. The explicit description of all irreducible isocrystals given in the introduction follows from Section \ref{sect:representatives}, Section \ref{sect:altdescofirred}, and Section \ref{sect:unramexts}. 

This explicit construction also extends to all indecomposable isocrystals (and thus gives explicit models of all isocrystals), following the lifting procedure of Section \ref{sect:representatives}.

The description of Kottwitz \cite[Lem.\ 11.2]{Kott} of the Hasse invariants of the endomorphism rings of simple objects of $\Isoc(\bF_q)$ is a special case of Corollary \ref{cor:descofendoinbrgroup}.

\subsection{Base Change For Isocrystals}\label{sect:basechangeisoc}

Recall that if $k$ is algebraically closed, $\Isoc(k)$ is semisimple, with simple objects $E_{\lambda}$ parametrised by $\lambda \in \bQ$ \cite{Manin}. Here writing $\lambda = \frac{a}{b} \in \bQ$, with $(a,b) = 1$ and $b \geq 1$, $E_{\lambda}$ is the isocrystal $(K(k)^{\oplus b}, F_{\lambda})$, where
\[
F_{\lambda}(\mathbf{v}) = R_b(p^a)\varphi(\mathbf{v}) = \begin{pmatrix}
    0 & 0 &\cdots & 0 & p^a \\
    1 & 0 & \cdots & 0 & 0 \\
    \vdots & \vdots & \ddots & \ddots & \vdots \\
    0 & 0 & \cdots & 0 & 0 \\
    0 & 0 & \cdots & 1 & 0
\end{pmatrix} \varphi(\mathbf{v}).
\]

For any characteristic $p$ extension $k' / k$ of perfect fields, there is a base change functor
\[
K(k') \otimes_{K(k)} (-) \colon \Isoc(k) \rightarrow \Isoc(k'), \qquad (X,F) \mapsto (K(k') \otimes_{K(k)} X, \varphi' \otimes F).
\]
For $k$ any perfect field, and $X \in \Isoc(k)$, one can use this to define the slopes of $X$ via base change to $\overline{k}$. Explicitly, if 
\[
K(\overline{k}) \otimes_{K(k)} X \cong \bigoplus_i E_{\lambda_i}^{\oplus m_i},
\]
then the slopes of $X$ are $\lambda_i$, where $\lambda_i$ is considered with multiplicity $m_i \cdot \dim_{K(\overline{k})}(E_{\lambda_i})$. From this, one obtains the Newton polygon, by ordering the slopes $\lambda_1 \leq \cdots \leq \lambda_k$ that appear with non-zero multiplicity, and letting $\text{NP}(X)$ to be the polygon starting at $(0,0)$ and consecutively starting from $\lambda_1$ adding in for each $\lambda_i$ a side of gradient $\lambda_i$ and horizontal projection length $m_i \dim_{K(\overline{k})}(E_{\lambda_i})$.

Now let us consider the case when $k = \bF_q$ is finite, where $q = p^r$ for some $r \geq 1$. As described in the introduction, in this case $K(\bF_q) = \bQ_q$ is the unique unramified extension of $\bQ_p$ of degree $r$, and $\varphi_q$ is the lift of the Frobenius automorphism of $\Gal(\bF_q / \bF_p)$. For $X \in \Isoc(\bF_q)$, one can show that the Newton polygon of $X$ is the same up to rescaling as the Newton polygon of the characteristic polynomial of the $\bQ_q$-linear endomorphism $F_X^r \colon X \rightarrow X$ \cite[Chap.\ IV, \S 6]{Dem}.

It is natural to ask how isocrystals behave under base change with respect to our explicit identification of the indecomposable objects of $\Isoc(\bF_q)$ of Section \ref{sect:basechange}
\[
(\Gal_{\bQ_p} \!\! \backslash \overline{\bQ_p}^{\times}) \times \bZ_{\geq 1} \xrightarrow{\sim} \Ind(\Isoc(\bF_q)), \qquad (\alpha, t) \mapsto X_{\alpha, t, q},
\]
where we write $X_{\alpha, t, q} \in \Isoc(\bF_q)$ for the indecomposable isocrystal corresponding to $(\alpha, t)$. To simplify notation, when $t = 1$ and $X_{\alpha, t, q}$ is simple, we write $X_{\alpha, q} \coloneqq X_{\alpha, 1, q}$.

\begin{eg}\label{eg:basechangeofsimplesfromFp}
Let $\bF \coloneqq \overline{\bF_p}$, which has $K(\bF) = \breve{\bQ}_p$ the completion of the maximal unramified extension of $\bQ_p$. From their explicit definition above, the simple objects $E_{\lambda} \in \Isoc(\bF)$ are all defined over $\bF_p$, as the definition matrix has values in $\bQ_p$, and so let us consider $E_{\lambda} \in \Isoc(\bF_p)$. In our notation, $E_{\lambda}$ is the isocrystal $X_{\beta, p} \in \Isoc(\bF_p)$, where $\beta$ denotes any fixed root of $x^b - p^a$ for $\lambda = a/b$, $b \geq 1$, $(a,b) = 1$. As this is irreducible over $\bQ_p$, all such roots are $\Gal_{\bQ_p}$-conjugate.

What is $\bQ_q \otimes_{\bQ_p} E_{\lambda} \in \Isoc(\bF_q)$? From Proposition \ref{prop:explicitdescofbasechange},
\[
\bQ_q \otimes_{\bQ_p} E_{\lambda} \cong \bQ_q \otimes_{\bQ_p} X_{\beta, p} \cong X_{\beta^r, q}^{N}, \qquad \text{ where } N = \frac{[\bQ_p(\beta) : \bQ_p(\beta^r)]}{m_{\bQ_q}^{\varphi_q}(\beta^r)},
\]
the formula for $N$ simplifying as the semilinear Schur indices for $\Isoc(\bF_p)$ are all trivial.

We can compute $N$ as follows. Write $r \lambda = a' / b'$ with $(a',b') = 1$, and $b' \geq 1$. Then $\beta^r$ is a root of $x^{b'} - p^{a'}$, so $x^{b'} - p^{a'}$ is the minimum polynomial of $\beta$ over $\bQ_q$ and thus the element of $\sP_{\bQ_q}$ corresponding to $X_{\beta^r, q}$. This has all coefficients in $\bQ_p$, hence from Corollary \ref{cor:unramSLSchur}
\[
m_{\bQ_q}^{\varphi_q}(\beta^r) = \frac{r}{\gcd(r, \nu_{\bQ_p(\beta)}(\beta))}
= \frac{r}{\gcd(r, a')}. 
\]
In particular, $m_{\bQ_q}^{\varphi_q}(\beta^r) > 1$ in general. On the other hand,
\[
[\bQ_p(\beta) : \bQ_p(\beta^r)] = \frac{[\bQ_p(\beta) : \bQ_p]}{[\bQ_p(\beta^r) : \bQ_p]} = \frac{b}{b'}.
\]
How are these related? Well, $r a b' = a' b$, so
\[
b \gcd(r, a') = \gcd(br, ba') = \gcd(br, rab') = r \gcd(b,ab') = r b',
\]
where $\gcd(b,ab') = b'$ because $(a,b) = 1$. In particular, $r / \gcd(r,a') = b / b'$, and thus $N = 1$.

In fact, this should come as no surprise: $\bQ_q \otimes_{\bQ_p} E_{\lambda}$ is irreducible, as its base change to $\Isoc(\bF)$ is, and therefore the only option is that $N = 1$.
\end{eg}

Now, we know from the Dieudonn\'{e}-Manin classification over $\bF$ that any $X \in \Isoc(\bF_q)$ becomes isomorphic to a direct sum of the simple modules $E_{\lambda}$ in $\Isoc(\bF)$. When is this isomorphism defined over a finite extension, in some $\Isoc(\bF_{q'})$ for $\bF_{q'} / \bF_q$? At first glance one might naturally guess that this is always true, as $\bF$ is the union of the finite extensions of $\bF_q$. The following gives a complete characterisation, which shows that the situation is essentially the opposite.

Let $\sC(\bF_q)$ denote the full subcategory of $\Isoc(\bF_q)$, whose objects are direct sums of the simple modules $\bQ_q \otimes_{\bQ_p} E_{\lambda}$ for $\lambda \in \bQ$. Let $\sD(\bF_q)$ denote the full subcategory of $\Isoc(\bF_q)$, whose objects are direct sums of the simple modules $X_{\alpha, 1, q}$ where $\alpha^m \in p^{\bZ}$ for some $m \geq 1$.

From Example \ref{eg:basechangeofsimplesfromFp}, the objects of $\sC(\bF_q)$ are the direct sums of the simple isocrystals $X_{\alpha, q}$ where $\alpha \in p^{\bQ}$, where we write 
\[
p^{\bQ} \coloneqq \{\alpha \in \overline{\bQ_p}^{\times} \mid \alpha^b = p^a \text{ for some } a, b \in \bZ \text{ with } (a,b) = 1 \text{ and } b \geq 1\}. 
\]
Indeed, Example \ref{eg:basechangeofsimplesfromFp} shows that $\bQ_q \otimes_{\bQ_p} E_{\lambda} \cong X_{\beta^r, q}$, where $\beta^r \in p^{\bQ}$. Conversely, if $\alpha \in p^{\bQ}$ satisfies $\alpha^{b'} = p^{a'}$ for some coprime $a', b' \in \bZ$ with $b' \geq 1$, then for $\lambda \coloneqq a' / r b'$, $\bQ_q \otimes_{\bQ_p} E_{\lambda} \cong X_{\alpha, q}$ by Example \ref{eg:basechangeofsimplesfromFp}.

Therefore we have that $\sC(\bF_q) \subset \sD(\bF_q)$, and $\sD(\bF_q)$ is the minor enlargement which adds in those objects $X_{\alpha, 1, q}$ for which $\alpha \zeta \in p^{\bQ}$ for some root of unity $\zeta$.

\begin{prop}
	Suppose that $X \in \Isoc(\bF_q)$. Then $\bQ_{q'} \otimes_{\bQ_q} X \in \sC(\bF_{q'})$ for some finite extension $\bF_{q'} / \bF_q$ if and only if $X \in \sD(\bF_q)$.
\end{prop}

\begin{proof}
	Because $\sD(\bF_q)$ and each $\sC(\bF_{q'})$ are closed under direct summands, it suffices to consider the case that $X$ is indecomposable, so $X = X_{\alpha, t, q}$. In this case, for any $\bF_{q'} / \bF_q$, $\bQ_{q'} \otimes_{\bQ_q} X$ is a direct sum of copies of $X_{\alpha^s, t, q'}$ by Proposition \ref{prop:explicitdescofbasechange}, where $s = r' / r$, writing $q = p^{r'}$. Now $X_{\alpha^s, t, q'} \in \sC(\bF_{q'})$ if and only if $t = 1$ and $\alpha^s \in p^{\bQ}$, which is the case if and only if $X \in \sD(\bF_q)$.
\end{proof}

\section{Semilinear Representations of Finite Cyclic Groups}\label{sect:cyclicgroups}

Suppose that $d \in \bZ_{\geq 1}$, and let $C_d = \langle z \rangle$ denote the cyclic group of order $d$. In this examine the consequences of the work of the previous chapters for the semilinear representation theory of $C_d$.

Let $\sigma$ be an automorphism of $L$ or order $r$ with $\sigma^d = 1$, so $r \mid d$. Then there is an exact sequence
\[
1 \rightarrow C_{d/r} \rightarrow C_d \xrightarrow{x \mapsto \sigma} \Gal(L / K) \rightarrow 1,
\]
where $K = L^{\sigma}$, and we can consider $\Rep_L^{\sigma}(C_d)$. The quotient map 
\[
\bZ \twoheadrightarrow C_d, \qquad 1 \mapsto z
\]
sends $r \bZ$ to the subgroup $C_{d/r}$, and inflation defines fully faithful exact inclusions
\[
\Rep_L(C_{d/r}) \hookrightarrow \Rep_L(r \bZ), \qquad \Rep_L^{\sigma}(C_d) \hookrightarrow \Rep_L^{\sigma}(\bZ)
\]
which preserve irreducibility and indecomposability. These are compatible with the induction and restriction functors between the respective linear and semilinear categories, and therefore the problem of computing semilinear Schur indices for $C_d$ is a special case of that for $\bZ$. 

More precisely, $\Ind_L(C_{d/r}) \subset \Ind_L(r \bZ)$, and under the identification $\Ind_L(r \bZ) \leftrightarrow (\Gamma \backslash \sP_L) \times \bZ_{\geq 1}$ of Section \ref{sect:semilinearbackground}, $\Ind_L(C_{d/r})$ corresponds to the subset $\Gamma \backslash \sP_{L, {d/r}}$ where
\[
\sP_{L, {d/r}} \coloneqq\{ (p(x),t) \in \sP_L \times \bZ_{\geq 1} \mid p(x)^t \mid (x^{d/r} - 1)\} \subset \sP_L \times \bZ_{\geq 1}.
\]
In particular, $\Rep_L^{\sigma}(C_d)$ is determined by $m_L^{\sigma}(p(x))$ for $p(x) \in \sP_L$ with $p(x) \mid x^{d/r} - 1$, by Theorem \ref{thm:mainthm}. The extra constraints this condition imposes on the $m_L^{\sigma}(p(x))$ are the following.

\begin{cor}\label{cor:mpxcyclicgroups}
	Suppose that $m \geq 1$, and $p(x) \in \sP_L$ with $p(x) \mid x^m - 1$. Let $\ord(\alpha)$ denote the order of $\alpha$, a $m$th root of $1$, so $\ord(\alpha) \mid m$. Then the prime factors of $m_L^{\sigma}(p(x))$ are a subset of the prime factors of $\ord(\alpha)$. In particular, if $\ord(\alpha)$ is coprime to $[L(\alpha) : K(\alpha)]$, then $m_L^{\sigma}(p(x)) = 1$. 
	
	If furthermore the period equals the index in $\Br(K(\alpha))$, then additionally $m_L^{\sigma}(p(x)) \mid \ord(\alpha)$.
\end{cor}

\begin{proof}
	We have that $m_L^{\sigma}(p(x)) = m_{L(\alpha)}^{\tau_{\alpha}}(x - \alpha)$, from which we may apply  \cite[Cor.\ 15.4]{TAYSchur} and \cite[Rem.\ 15.5]{TAYSchur} because $x - \alpha$ corresponds to a representation of dimension $1$. This gives the relations between $m_L^{\sigma}(p(x))$ and $\ord(\alpha)$. As is true for any $p(x) \in \sP_L$, $m_L^{\sigma}(p(x)) \mid [L(\alpha) : K(\alpha)]$ by \cite[Prop.\ 9.3(6)]{TAYSchur}, and therefore $m_L^{\sigma}(p(x)) = 1$ if $\ord(\alpha)$ is coprime to $[L(\alpha) : K(\alpha)]$.
\end{proof}

\begin{remark}
	When the period equals the index in $\Br(K(\alpha))$, one can deduce that $m_L^{\sigma}(p(x)) \mid \ord(\alpha)$ directly using Theorem \ref{thm:mainthm}, as $\alpha^{\ord(\alpha)} = 1$ in $K(\alpha)^{\times} / N_{L(\alpha) / K(\alpha)}(L(\alpha)^{\times})$.
\end{remark}

The conditions of Corollary \ref{cor:mpxcyclicgroups} put strong constraints on the $m(p(x))$, but for specific extensions one can often say more.

\begin{eg}
	Suppose that $L / K$ is $\bC / \bR$, with $\sigma = \overline{(-)}$. Then $r = 2$, $d$ is even, and $\Rep_{\bC}(C_{d/2})$ has $\Irr_{\bC}(C_{d/2})$ parametrised by the $d / 2$th roots $\zeta$ of $1$ in $\bC$, corresponding to the irreducible polynomials $p(x) = x - \zeta$. From Section \ref{eg:CoverR}, we have that $m_{\bC}^{\sigma}(x - \zeta) = 1$ when $\zeta \neq -1$, and $m_{\bC}^{\sigma}(x + 1) = 2$. This is stronger than what one can deduce solely using Corollary \ref{cor:mpxcyclicgroups}.
\end{eg}
	
	Things are also more tractable for characteristic zero local fields.
	
	\begin{prop}
		Suppose that $L$ is a finite extension of $\bQ_p$ for some prime $p$, $L / K$ has ramification index $e = e_t p^a$, where $e_t$ is coprime to $p$. Suppose that $\zeta \in \mu_{p'}(K)$ is a root of $1$ in $K$ of order coprime to $p$. Then $\zeta \in N_{L/K}(L^{\times})$ if and only if $\zeta \in \mu_{p'}(K)^e$, if and only if $\zeta \in \mu_{p'}(K)^{e_t}$.
        
        In particular, $m_L^{\sigma}(x - \zeta)$ is the order of $\zeta$ in $\mu_{p'}(K) / \mu_{p'}(K)^e$.
	\end{prop}
	
	\begin{proof}
		The Artin map $\phi_{L/K} \colon K^{\times} \twoheadrightarrow \Gal(L/K)$ is a surjection, with kernel $N_{L/K}(L^{\times})$, which restricts to a surjection $\phi_{L/K} \colon \OO_K^{\times} \twoheadrightarrow I(L/K)$ to the inertia subgroup of $\Gal(L/K)$ of order $e$. Let $P(L/K) \subset I(L/K)$ denote the wild inertia subgroup, the unique Sylow $p$-subgroup of $I(L/K)$, and $I_t(L/K) \coloneqq I(L/K) / P(L/K)$ the tame inertia, of order $e_t$. The unit group $\OO_K^{\times} = \mu_{p'}(K) \times (1 + \fm_K)$ is a direct product, where $(1 + \fm_K)$ is pro-$p$, and thus composition of $\phi_{L/K}$ with the quotient map $I(L/K) \rightarrow I_t(L/K)$ induces a surjection
		\[
		\overline{\phi_{L/K}} \colon \mu_{p'}(K) \twoheadrightarrow I_t(L/K)
		\]
		with $\ker(\overline{\phi_{L/K}}) = \mu_{p'}(K) \cap \ker(\phi_{L/K})$. Because $\mu_{p'}(K)$ is cyclic, this subgroup is thus the unique subgroup of index $e_t$, $\mu_{p'}(K)^{e_t}$. As $\mu_{p'}(K)$ has order coprime to $p$, $\mu_{p'}(K)^{e} = \mu_{p'}(K)^{e_t}$.
		
		The final statement follows from Corollary \ref{cor:schurindexisorder} and Remark \ref{rem:localandglobBrauerGroups}.
	\end{proof}

\section{Algorithms}\label{sect:algorithms}

In this final section we briefly summarise the various algorithms that we have given throughout the text. We include magma implementations of them over number fields in the repository \cite{GITHUB}.

\begin{enumerate}
	\item Finding the generalised Jordan normal form of a matrix, Section \ref{sect:GJNF}. This is not original, but needed for our other algorithms.
	\item Testing if two matrices are $\sigma$-conjugate, Corollary \ref{cor:linsimimpliessim}, reducing to testing if two matrices are conjugate in the usual sense, which is standard, for example using algorithm (1).
	\item Given representatives for the semilinear normal form for orbits with $t = 1$, lift these to representatives for all $t \geq 1$. Section \ref{sect:representatives}.
	\item If $L$ is a local or global field, compute the semilinear Schur index of any $p(x) \in \sP_L$. This uses Theorem \ref{thm:mainthm}, as the order of elements in the quotient of norm groups that appears can be computed over these fields.
	\item For $X \in \GL_n(L)$, without having fixed precise representatives $A_{\OO}$, find the multiplicities that each $A_{\OO}$ appears in the semilinear normal form of $A$. Section \ref{sect:matrixsemilinearpov}. This relies on a knowledge of the semilinear Schur indices, and so is possible over local and global fields using algorithm (4).
	\item Find explicit representatives $A_{\OO}$ for the semilinear normal form, over unramified extensions of local fields, finite fields, and global fields $L$ where $L/K$ has prime degree. Section \ref{sect:representatives}.
\end{enumerate}

\bibliography{biblio}{}
\bibliographystyle{plain}
\end{document}